\documentclass[11pt]{article}
\usepackage{amssymb,amsmath,latexsym}
\usepackage[english]{babel}
\usepackage[T1]{fontenc}
\usepackage{url}
\usepackage{color}
\usepackage{ulem}
\usepackage[dvips]{graphics}

\definecolor{myblue}{rgb}{0,0,1}

\newcommand{\boldeta}{\boldsymbol{\eta}}
\newcommand{\bgamma}{\boldsymbol{\gamma}}
\newcommand{\balpha}{\boldsymbol{\alpha}}
\newcommand{\bbeta}{\boldsymbol{\beta}}

\newcommand{\bepsilon}{\boldsymbol{\epsilon}}

\newcommand{\bzeta}{\boldsymbol{\zeta}}
\newcommand{\blambda}{\boldsymbol{\lambda}}

\newcommand{\btheta}{\boldsymbol{\theta}}
\newcommand{\bvartheta}{\boldsymbol{\vartheta}}

\newcommand{\bPsi}{\boldsymbol{\Psi}}

\newcommand{\bzero}{\boldsymbol{0}}
\newcommand{\bnu}{\boldsymbol{\nu}}
\newcommand{\bmu}{\boldsymbol{\mu}}

\newcommand{\bSigma}{\boldsymbol{\Sigma}}
\newcommand{\bsigma}{\boldsymbol{\sigma}}
\newcommand{\bomega}{\boldsymbol{\omega}}
\newcommand{\bOmega}{\boldsymbol{\Omega}}

\newcommand{\bb}{\boldsymbol{b}}

\newcommand{\bA}{\boldsymbol{A}}

\newcommand{\bC}{\boldsymbol{C}}

\newcommand{\bD}{\boldsymbol{D}}

\newcommand{\bI}{\boldsymbol{I}}
\newcommand{\bJ}{\boldsymbol{J}}

\newcommand{\bcalI}{\boldsymbol{\cal I}}
\newcommand{\bcalJ}{\boldsymbol{\cal J}}
\newcommand{\bcalK}{\boldsymbol{\cal K}}
\newcommand{\bcalL}{\boldsymbol{\cal L}}
\newcommand{\bK}{\boldsymbol{K}}

\newcommand{\bx}{\boldsymbol{x}}

\newcommand{\by}{\boldsymbol{y}}

\newcommand{\bz}{\boldsymbol{z}}
\newtheorem{theo}{Theorem}
\newtheorem{lem}{Lemma}
\newtheorem{prop}{Proposition}

\numberwithin{equation}{section}
\numberwithin{theo}{section}
\numberwithin{lem}{section}
\numberwithin{cor}{section}
\numberwithin{exemple}{section}
\numberwithin{rem}{section}
\numberwithin{prop}{section}
\small\normalsize

\def\zak{\null\hfill{$\Box$}\par\vspace*{0.2cm}}

\def \argmin{\mathop{\hbox{\rm arg min}}}

\def\1g{1\hskip -3pt \mbox{l}}

\title{Goodness-of-fit tests for log and exponential GARCH models}
\author{
{\sc Christian Francq\footnote{CREST and University Lille 3
(EQUIPPE), BP 60149, 59653 Villeneuve d'Ascq cedex, France. E-Mail:
christian.francq@univ-lille3.fr}, Olivier Wintenberger\footnote{Universities of Paris 6 and Copenhagen, LSTA 4 Place Jussieu, 75005 Paris, France. E-Mail:
olivier.wintenberger@upmc.fr} and Jean-Michel
Zakoïan\footnote{Corresponding author: Jean-Michel Zakoïan, EQUIPPE
(University Lille 3) and CREST, 15 boulevard Gabriel Péri, 92245
Malakoff Cedex, France. E-mail: zakoian@ensae.fr, Phone number:
33.1.41.17.77.25.
 } }}
 \date{}
\begin{document}




\maketitle

\begin{abstract}
This paper studies goodness of fit tests and specification tests for an extension of the Log-GARCH model which is both asymmetric and stable by scaling.
A Lagrange-Multiplier  test is derived for testing
the extended Log-GARCH against more general formulations taking the form of combinations of
Log-GARCH and  Exponential GARCH (EGARCH). The null assumption of  an
EGARCH is also tested. Portmanteau goodness-of-fit tests are
developed for the extended Log-GARCH. An application to real financial data is
proposed.
% \PACS{PACS code1 \and PACS code2 \and more}\\
\end{abstract}
It is now widely accepted that, to model the dynamics of daily financial returns, volatility models have to incorporate the so-called leverage effect.\footnote{This effect, typically observed on most
stock returns series, means that negative returns have more impact on the volatility than positive returns of the same magnitude.}
Among the various asymmetric GARCH processes introduced in the econometric literature, E(xponential)GARCH and Log-GARCH models share the
property of specifying the dynamics of the log-volatility, rather than the volatility,
as a linear combination of past variables.
%values and past values of the positive and negative parts of the
One advantage of such specifications is to avoid positivity constraints on the parameters, which complicate statistical inference of standard GARCH formulations.
A class of (asymmetric) Log-GARCH(p,q) models was recently
studied by Francq, Wintenberger and Zako\"ian  (2013) (FWZ). In
this class, originally introduced by
Geweke (1986), Pantula (1986) and Milhøj (1987) (see Sucarrat, Grønneberg and Escribano (2015) for a more recent reference), the dynamics is
defined by
\begin{equation}\label{logGARCH}
\left\{\begin{array}{lll}
\epsilon_t&=&\sigma_t\eta_t,\\
\log \sigma_t^2&=&\omega+\sum_{i=1}^q\left(\alpha_{i+}1_{\{\epsilon_{t-i}>0\}}+\alpha_{i-}1_{\{\epsilon_{t-i}<0\}}\right)\log\epsilon_{t-i}^2\\
&&+\sum_{j=1}^p\beta_j\log \sigma_{t-j}^2
\end{array}\right.
\end{equation} where $\sigma_t>0$ and $(\eta_t)$ is a
sequence of independent and identically distributed (iid) variables
such that %$E\eta_1=0$ and
$E\eta_1^2=1$.

%The usual symmetric Log-GARCH corresponds to the case
%$\balpha_{+}=\balpha_{-}$, with
%$\balpha_{+}=(\alpha_{1+},\dots,\alpha_{q+})$ and
%$\balpha_{-}=(\alpha_{1-},\dots,\alpha_{q-})$.  An advantage of
%modeling the log-volatility rather than the volatility is that the
%vector $\btheta=(\omega,\balpha_+,\balpha_-,\bbeta)$ with
%$\bbeta=(\beta_{1},\dots,\beta_{p})$ is not a priori subject to
%positivity constraints\footnote{However, some desirable properties
%may determine the sign of  coefficients. For instance, the present
%volatility is generally thought of as an increasing function of its
%the past values, which entails $\beta_j>0.$ The difference with
%standard GARCH models is that such constraints are not required for
%the existence of the process and, thus, do not complicate the
%estimation procedures.}.

One drawback of this model is that it is generally not stable by scaling.
Indeed, if $(\epsilon_t)$ is a solution of Model (\ref{logGARCH}),
the process $(\epsilon_t^*)$ defined by $\epsilon_t^*=c\epsilon_t$
with $c>0$ satisfies $\epsilon_t^*=\sigma_t^*\eta_t$ with $\sigma_t^{*2}=\omega_{t-1}^*+\sum_{i=1}^q\left(\alpha_{i+}1_{\{\epsilon_{t-i}^*>0\}}+\alpha_{i-}1_{\{\epsilon_{t-i}^*<0\}}\right)
\log\epsilon_{t-i}^{*2}
+\sum_{j=1}^p\beta_j\log \sigma_{t-j}^{*2}$
%
%$$%\label{logGARCHscale}
%\left\{\begin{array}{lll}
%\epsilon_t^*&=&\sigma_t^*\eta_t,\quad \\
%\log \sigma_t^{*2}&=&\omega_{t-1}^*+\sum_{i=1}^q\left(\alpha_{i+}1_{\{\epsilon_{t-i}^*>0\}}+\alpha_{i-}1_{\{\epsilon_{t-i}^*<0\}}\right)\log\epsilon_{t-i}^{*2}\\
%&&+\sum_{j=1}^p\beta_j\log \sigma_{t-j}^{*2}
%\end{array},\right.
%$$
where
$$\omega_{t-1}^*=\log c^2\left(1-\sum_{j=1}^p\beta_j-\sum_{i=1}^q\left(\alpha_{i+}1_{\{\epsilon_{t-i}^*>0\}}+\alpha_{i-}1_{\{\epsilon_{t-i}^*<0\}}\right)\right)$$
is not constant (except in the symmetric case where
$\alpha_{i+}=\alpha_{i-}$ for all $i$). It is important that a volatility model be stable
by scaling.\footnote{Indeed, as remarked by a referee,
a practitioner is essentially faced by three choices:
(a) leave returns untransformed, i.e. set $c = 1$, (b) express returns in
terms of percentages, i.e. set $c = 100$, or (c) express returns in
terms of basis points, i.e. set $c = 10,000$. Clearly, it is desirable
that the dynamics of the volatility model be not affected by the choice of $c$.}
The standard log-GARCH has the stability by scaling property, but is not able to capture the leverage effect.

In this paper, we will consider an extension of Model (\ref{logGARCH}) which is both stable by scaling and asymmetric.
Our main foci concern specification tests of this model and the comparison with the EGARCH model.
%In this paper, our main focus is on specification tests for Log-GARCH and EGARCH dynamics. %We also wish to develop goodness-of-fit tests for the extended Log-GARCH model.
%A natural alternative of the Log-GARCH is the EGARCH model
The latter formulation, introduced by Nelson (1991), appears as a widely used competitor of the Log-GARCH in applications.  As we will see, the two models display very similar properties
and their volatility dynamics may coincide. However, %even if
%Log-GARCH and EGARCH models have apparent similarities and certain common properties,
the Log-GARCH and EGARCH models are not equivalent from a statistical point of view.
In particular, it is obvious to invert the Log-GARCH model, {\it i.e.} to express the volatility as an explicit function of the past returns,
whereas the EGARCH(1,1) is invertible only under strong restrictions on the parameters. This is a major drawback for the statistical inference of the second specification, see Wintenberger (2013) and FWZ.
%So one could prefer the first specification, as the Log-GARCH model can produce the same volatility process as the EGARCH %(see Lemma~\ref{EGARCHisLogGARCH})
%and the inference is much easier with the Log-GARCH.
However, the two models are not compatible for a same series and one has to discuss if one specification is more likely to fit the data at hand than the other.
It is therefore of interest to develop testing procedures for one specification against the other. This constitutes the main aim of the present paper.

The remainder of the paper is organized as follows.
Section~\ref{sec1} introduces the extended Log-GARCH model and discusses its similarities with the EGARCH.
It also provides strict stationarity conditions.
Section~\ref{sec2} studies %the existence of strictly stationary solutions to Model (\ref{logGARCH}) and
the asymptotic properties of the quasi-maximum likelihood
(QML) estimator. Section~\ref{sec4} considers testing the null
assumption of a Log-GARCH against more general formulations
including the EGARCH. Section~\ref{sec5}
considers the reverse problem, in which the null assumption is the
EGARCH model. In Section~\ref{sec6}, Portmanteau goodness-of-fit
tests are developed for the Log-GARCH. Section~\ref{sec7} compares the Log-GARCH and EGARCH models for series of exchange rates.
%some numerical applications on simulated and real data. %Section
%\ref{sec8} concludes.

\section{Extended Log-GARCH model}
\label{sec1}
Consider the {\it Asymmetric and stable by Scaling}
Log-GARCH (AS-Log-GARCH)  model of order  ($p,q$), defined by
\begin{equation}\label{logGARCHSTABLE}
\left\{\begin{array}{lll}
\epsilon_t&=&\sigma_t\eta_t,\quad \\
\log \sigma_t^2&=&\omega+\sum_{i=1}^q\omega_{i-}1_{\{\epsilon_{t-i}<0\}}+\sum_{j=1}^p\beta_j\log \sigma_{t-j}^2\\
&&+\sum_{i=1}^q\left(\alpha_{i+}1_{\{\epsilon_{t-i}>0\}}+\alpha_{i-}1_{\{\epsilon_{t-i}<0\}}\right)\log\epsilon_{t-i}^2,
\end{array}\right.
\end{equation}
where $\omega$ and the components of the vectors $\bomega_{-}=(\omega_{1-},\dots,\omega_{q-})'$, $\balpha_{+}=(\alpha_{1+},\dots,\alpha_{q+})'$,
$\balpha_{-}=(\alpha_{1-},\dots,\alpha_{q-})'$, and
$\bbeta=(\beta_{1},\dots,\beta_{p})'$ are
real coefficients, which are not {\it a priori} subject to
positivity constraints, under the same assumptions on $(\eta_t)$ as
in Model (\ref{logGARCH}).
The main features of the asymmetric Log-GARCH($p,q$) model -
volatility which is not bounded below, persistence of small values,
power-aggregation - continue to hold in
this extended version. We refer the reader to FWZ  for details.
Contrary to Model (\ref{logGARCH}), the extended formulation
(\ref{logGARCHSTABLE}) is stable by scaling. Moreover, this model
leads to a different interpretation of the usual leverage effect.

\subsection{News Impact Curves}
\begin{center}
\begin{figure}%[h]
\vspace*{7.cm} \protect
\includegraphics{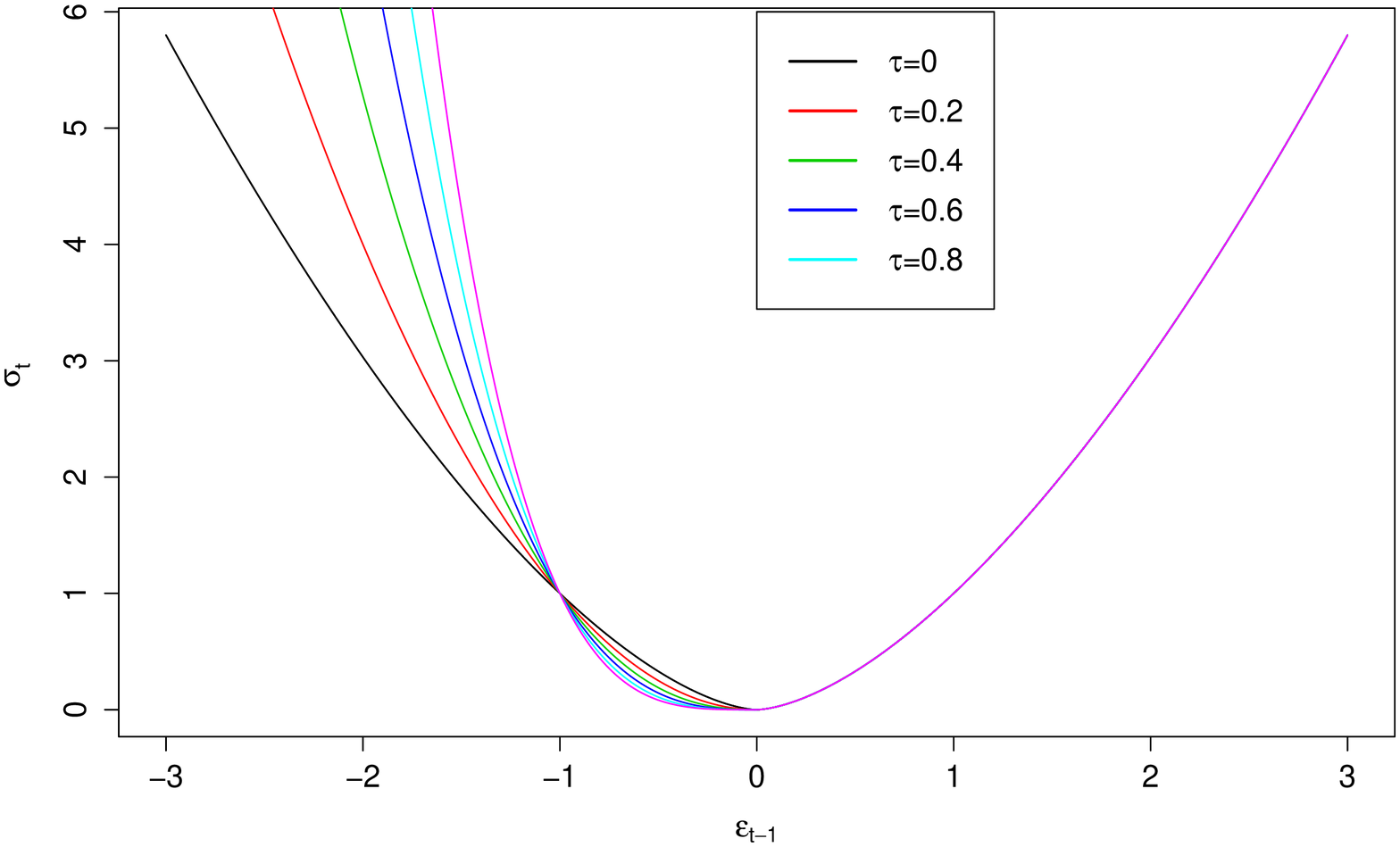}
\includegraphics{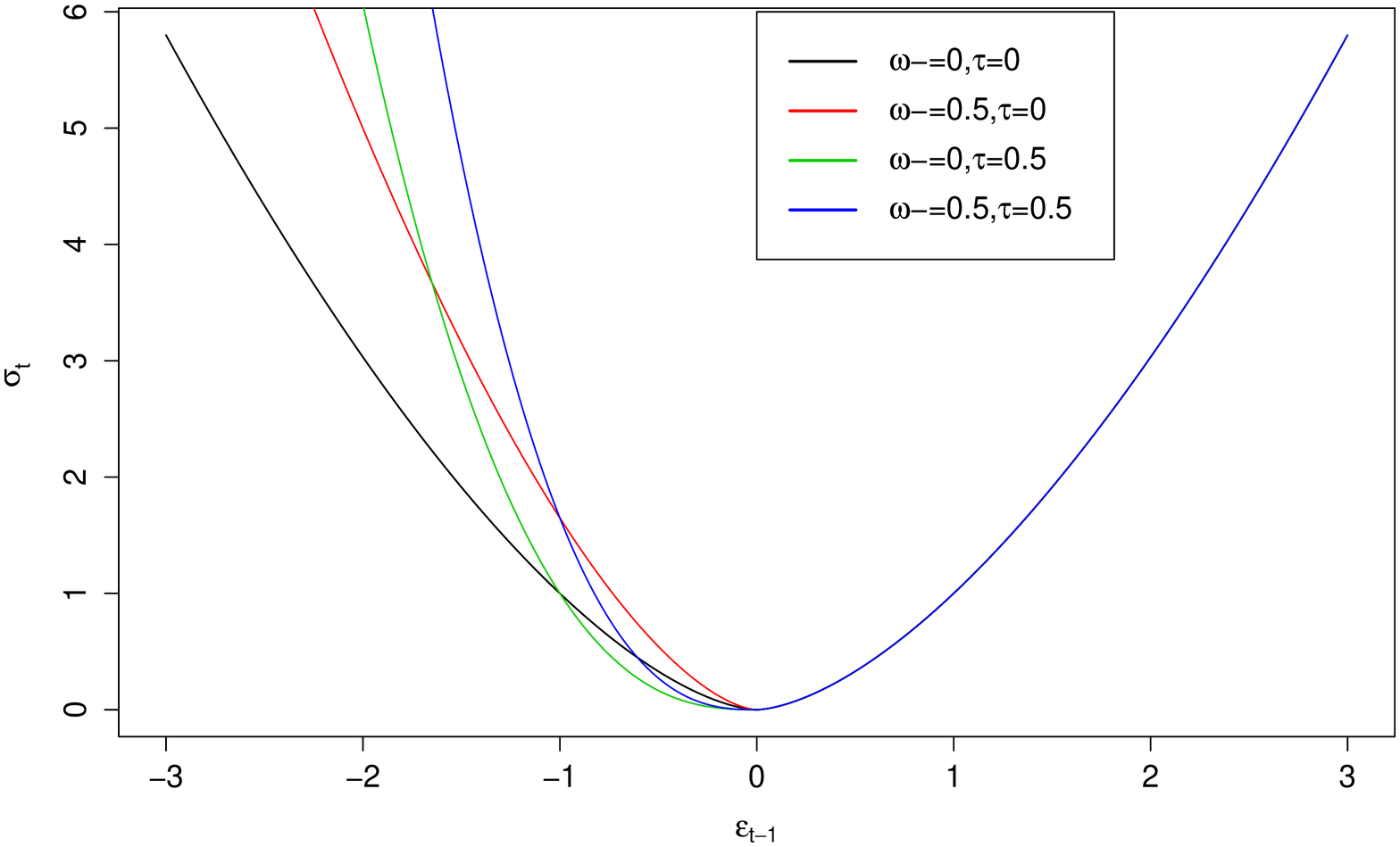}
\includegraphics{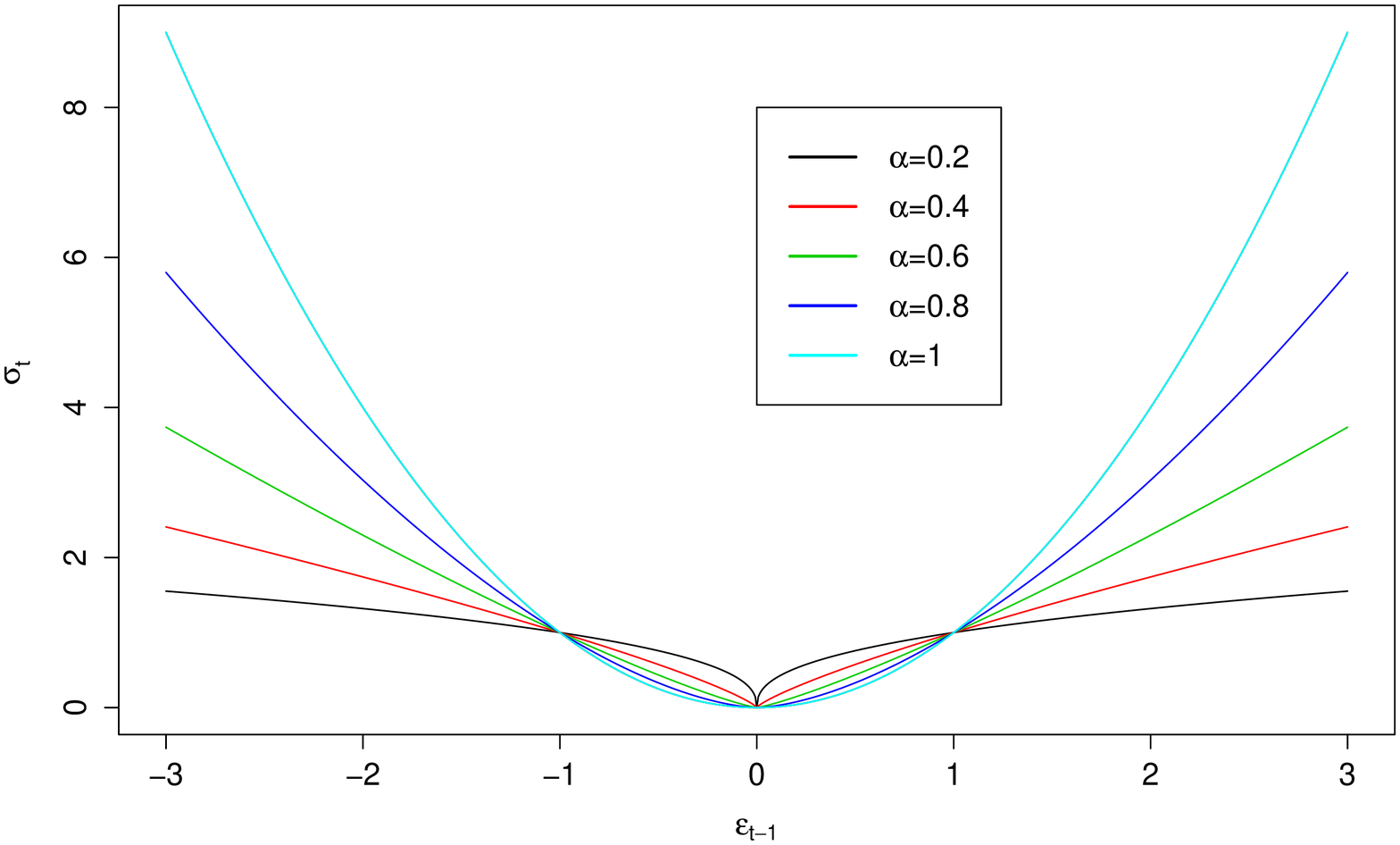}
\includegraphics{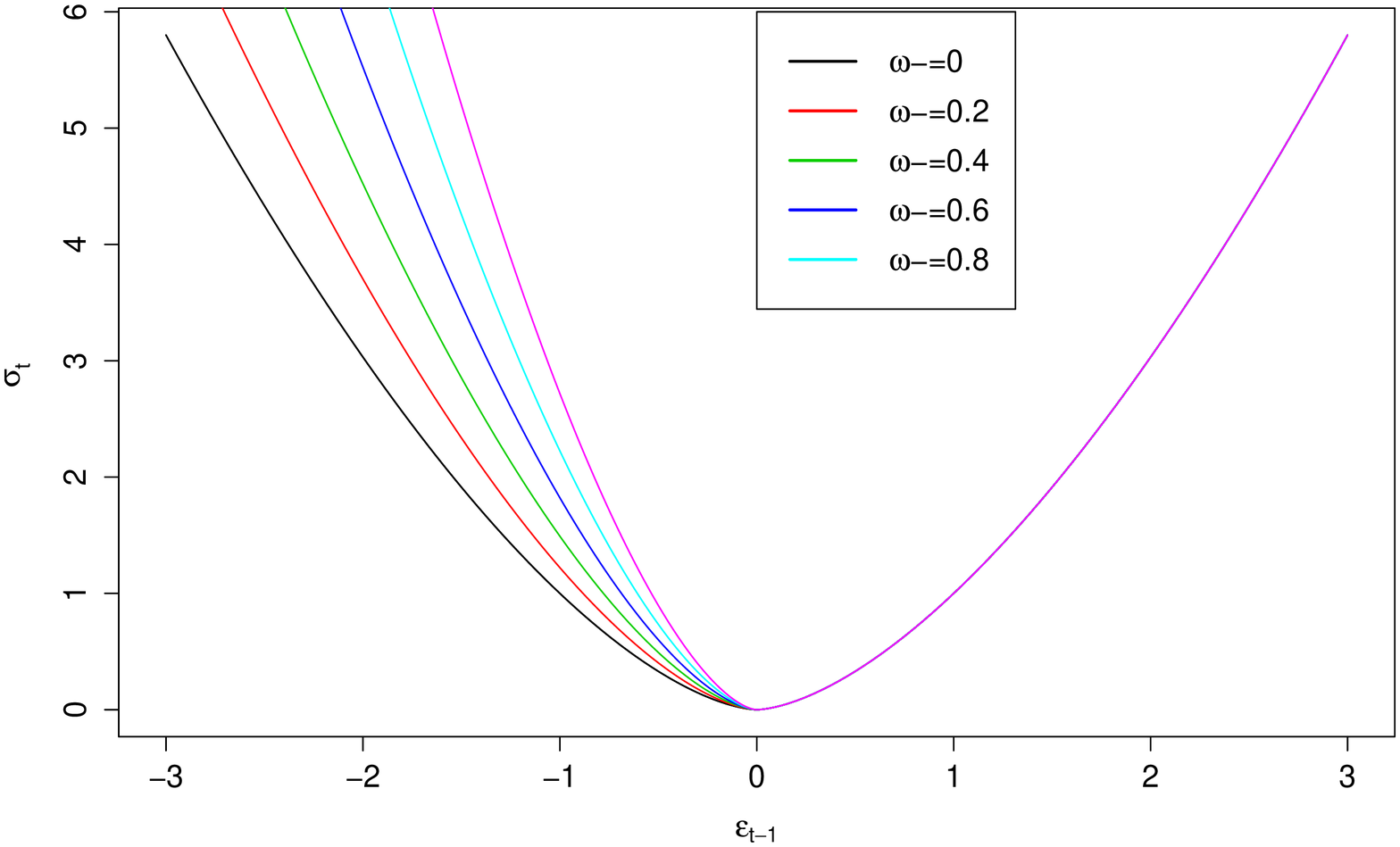}

\vspace*{1cm}
\caption{\label{fig1} {\small News Impact Curves: $\sigma_t$ as a function of $\epsilon_{t-1}$ in (\ref{tau}). The parameter $\omega$ is set to 0.
%The top left graph is obtained for $\tau=0$, the top right graph for $\alpha=0.8$, and the bottom left graph for $\omega_-=0$.
The top graphs are obtained for $\tau=0$, the left graphs for $\omega_-=0$, the right graphs and the bottom left graph for $\alpha=0.8$.
}}
\end{figure}
\end{center}

Compared to model \eqref{logGARCH}, the AS-Log-GARCH model (\ref{logGARCHSTABLE}) contains additional asymmetry parameters.
Through the introduction of the coefficients $\omega_{i-}$, Model
(\ref{logGARCHSTABLE})
 allows for an asymmetric impact of the past positive
and negative returns on the log-volatility which does not depend on
their magnitudes. For instance, consider the AS-Log-ARCH(1)
model with $\alpha_{1+}= \alpha_{1-}=\alpha.$ We have
$$\sigma_t^2= e^{\omega+ \omega_{1-}1_{\{\epsilon_{t-1}<0\}}}(\epsilon_{t-1}^2)^{\alpha}.$$
If $\omega_{1-}>0$, a decrease of the price, whatever its amplitude,
will increase the volatility by a scaling factor $e^{ \omega_{1-}}.$
In the limit case where $\alpha=0$, the volatility takes only two values depending only on the sign (not the size) of the past return.
Now we turn to the second leverage effect.
If $\alpha_{1+}=\alpha$ and $\alpha_{1-}=\alpha+\tau$ with
$\tau>0$, we have
\begin{equation}\label{tau}
\sigma_t^2= e^{\omega+ \omega_{1-}1_{\{\epsilon_{t-1}<0\}}}(\epsilon_{t-1}^2)^{\alpha}
(\epsilon_{t-1}^2)^{\tau
1_{\{\epsilon_{t-1}<0\}}}.\end{equation} The effect of
a  large negative return ($\epsilon_{t-1}<-1$) is an increase of
volatility, but the effect may be reversed for very small returns.
For small but not too small returns, this effect is balanced by the
presence of the scaling factor $e^{\omega_{1-}}$.
To summarize, the AS-Log-GARCH is in fact capable of
detecting two types of leverage: one type where the leverage effect depends on the magnitude
of negative return, and one type in which it does not.
The so-called News Impact Curves, displaying $\sigma_t$ as a function of $\epsilon_{t-1}$, are provided in Figure \ref{fig1}.
%\begin{rem}\label{ident?}
%One could think of a different parametrization
%$$
%\left\{\begin{array}{lll}
%\epsilon_t&=&\sigma_t\eta_t,\quad \\
%\log \sigma_t^2&=&\sum_{i=1}^q\left((\omega_{i+}'1_{\{\epsilon_{t-i}> 0\}}+\omega_{i-}'1_{\{\epsilon_{t-i}<0\}})+(\alpha_{i+}1_{\{\epsilon_{t-i}>0\}}+\alpha_{i-}1_{\{\epsilon_{t-i}<0\}}\log\epsilon_{t-i}^2)\right)\\
%&&+\sum_{j=1}^p\beta_j\log \sigma_{t-j}^2,
%\end{array}\right.
%$$
%where $\omega_{i+}', \omega_{i-}',\alpha_{i+}, \alpha_{i-}, \beta_j$ are
%real coefficients. However, this parametrization is not identifiable as it coincides with the model \eqref{logGARCHSTABLE} for $\omega=\sum_{i=1}^q\omega_{i+}'$ and $\omega_{i-}=\omega_{i-}'- \omega_{i+}'$, $1\le i\le q$. We will see in the sequel that the parametrization \eqref{logGARCHSTABLE} is identifiable.\end{rem}

\subsection{Similarities with the EGARCH dynamics}
%In this paper, our main focus is on specification tests for Log-GARCH and EGARCH dynamics. %We also wish to develop goodness-of-fit tests for the extended Log-GARCH model.
%%A natural alternative of the Log-GARCH is the EGARCH model
%The latter model, introduced by Nelson (1991), appears as a competitor of the Log-GARCH in applications as the two models display at a first sight very similar dynamics.
The dynamics of the logarithm of the volatility of the EGARCH$(p,\ell)$ model is provided by the recursion
\begin{equation}
\label{egarch}
\log \sigma_t^2=\tilde\omega+\sum_{j=1}^p\tilde \beta_{j}\log \sigma_{t-j}^2
+ \sum_{k=1}^{\ell} \gamma_{k+}\tilde\eta^{+}_{t-k} + \gamma_{k-}\tilde\eta^{-}_{t-k},
\end{equation}
where the innovations $\tilde \eta_t$ are iid random variables such that $E\tilde \eta^2_1=1$, with the notation $x^+=\max\{x,0\}$ and $x^-=\max\{-x,0\}$.
If one substitutes $\log\sigma_{t-i}^2+\log\eta_{t-i}^2$ for $\log\epsilon_{t-i}^2$ in \eqref{logGARCHSTABLE}, the probabilistic structures of the two classes of models seem similar.
More precisely, we have the following result.
\begin{prop}
\label{EGARCHisLogGARCH}

(i) For any EGARCH process $\tilde{\epsilon}_t=\sigma_t\tilde{\eta}_t$ satisfying \eqref{egarch} with $Ee^{s_0|\tilde{\eta}_1}|<\infty$ for some $s_0>0$, there exists a AS-Log-GARCH process  $\epsilon_t=\sigma_t\eta_t$ satisfying \eqref{logGARCHSTABLE}, with the same volatility process $\sigma_t$ and $\eta_t$ measurable with respect to $\tilde{\eta}_t$. %for which the volatility processes of $(\epsilon_t)$ and $(\tilde{\epsilon}_t)$ are the same.

(ii) Conversely, there exist AS-Log-GARCH processes $\epsilon_t=\sigma_t\eta_t$ for which there is no EGARCH process $\tilde{\epsilon}_t=\sigma_t\tilde{\eta}_t$
with the same volatility process $\sigma_t$ and $\tilde{\eta}_t$ measurable with respect to $\eta_t$.
\end{prop}
{\bf Proof:} Let us prove (i). For simplicity of notation, we assume that $\tilde{\epsilon}_t=\sigma_t\tilde{\eta}_t$ follows the first order EGARCH$(1,1)$ model, and we drop the indexes $i$, $j$ and $k$. Let the Log-GARCH(1,1) process $\epsilon_t=\sigma_t\eta_t$ satisfying \eqref{logGARCHSTABLE} with  the parameters $\alpha:=\alpha+=\alpha_-\neq 0$, $\omega+\alpha c_+=\tilde\omega$, $\omega_{-}=-\alpha(c_--c_+)$,  and $\alpha+\beta=\tilde{\beta}$, and the noise   $\eta_t=e^{\frac{c+}{2}}e^{\frac{\gamma_+}{2\alpha}|\tilde{\eta}_t|}1_{\tilde{\eta}_t\geq 0}-e^{\frac{c-}{2}}e^{\frac{\gamma_-}{2\alpha}|\tilde{\eta}_t|}1_{\tilde{\eta}_t< 0},$
with constants $c_+$ and $c_-$ to be chosen later. The Log-GARCH volatility then satisfies
\begin{eqnarray*}\log\sigma_t^2&=&\omega+\omega_-1_{\eta_{t-1}<0}+\alpha\log\eta_{t-1}^2+(\alpha+\beta)\log\sigma_{t-1}^2\\
&=&\tilde{\omega}+(\gamma_+1_{\tilde{\eta}_{t-1}>0}+\gamma_-1_{\tilde{\eta}_{t-1}<0})|\tilde{\eta}_{t-1}|+\tilde{\beta}\log\sigma_{t-1}^2,
\end{eqnarray*}
which is the equation satisfied by the volatility of the EGARCH(1,1) model. It then suffices to choose $\alpha$ such that $\gamma_+/\alpha<s_0$ and $\gamma_-/\alpha<s_0$, and then  $c_+$ and $c_-$ such that $E\eta_t^2=1$.

Now we turn to (ii). Let $(\epsilon_t)$ denote any AS-Log-GARCH process satisfying \eqref{logGARCHSTABLE}, with $\alpha_{1+}\ne \alpha_{1-}$, and sufficiently general so that
the support of the law of $\log\sigma_{t-1}^2$ contains at least three different values. Also assume that
$\log\eta_{t-1}^2$ has a finite variance. We proceed by contradiction.
Suppose there exists an EGARCH process
satisfying $\tilde{\epsilon}_t=\sigma_t\tilde{\eta}_t$ with $\tilde{\eta}_t=f(\eta_t)$ for some measurable function $f$. We thus have
\begin{eqnarray*}\log\sigma_t^2&=&\omega+\omega_-1_{\eta_{t-1}<0}+
\left(\alpha_{1+}1_{\{\epsilon_{t-1}>0\}}+\alpha_{1-}1_{\{\epsilon_{t-1}<0\}}\right)\log\eta_{t-1}^2\\&&
+ \left(\alpha_{1+}1_{\{\epsilon_{t-1}>0\}}+\alpha_{1-}1_{\{\epsilon_{t-1}<0\}}\right)\log\sigma_{t-1}^2+ \sum_{j=1}^p\beta_j\log \sigma_{t-j}^2\\
&=&\tilde{\omega}+(\gamma_+1_{\tilde{\eta}_{t-1}>0}+\gamma_-1_{\tilde{\eta}_{t-1}<0})|\tilde{\eta}_{t-1}|+\tilde{\beta}\log\sigma_{t-1}^2\\&&
+\sum_{j=1}^p\tilde \beta_{j}\log \sigma_{t-j}^2
+ \sum_{k=2}^{\ell} \gamma_{k+}\tilde\eta^{+}_{t-k} + \gamma_{k-}\tilde\eta^{-}_{t-k},
\end{eqnarray*}
which entails
$$a(\eta_{t-1})=b_{t-2}+c(\eta_{t-1})\log\sigma_{t-1}^2$$
where $b_{t-2}$ denotes a variable belonging to $\sigma$-field ${\cal F}_{t-2}$ generated by the $\eta_{t-2-j}$  with $j\geq 0$.
%Because $(a_{t-1},c_{t-1})$ is independent from $(b_{t-2},d_{t-2})$, it follows that
We have
\begin{eqnarray*}0&=&\mbox{var}\{a(\eta_{t-1})-b_{t-2}-c(\eta_{t-1})\log\sigma_{t-1}^2|{\cal F}_{t-2}\}\\&=&
\mbox{var}\{a(\eta_{t-1})\}+\log^2\sigma_{t-1}^2\mbox{var}\{c(\eta_{t-1})\}-2\log\sigma_{t-1}^2\mbox{cov}\{a(\eta_{t-1}), c(\eta_{t-1})\},
\end{eqnarray*}
from which it follows that $\log\sigma_{t-1}^2$ takes at most two values.
This contradicts the above assumptions.
\zak
This proposition allows to complete the interpretation of the two types of
       leverage effects in the  AS-Log-GARCH. The coefficients $\omega_{0,i-}$ produce the leverage effect of the EGARCH volatility, i.e. an asymmetry depending on the amplitude of the innovations $\tilde{\eta}_{t-i}$. On the opposite, the EGARCH model cannot capture the asymmetric effect induced by the coefficients $\alpha_{0,i-}, \alpha_{0,i+}$ and the amplitude of the returns $\epsilon_{t-i}$. %This last leverage effect is also captured by asymmetric GARCH models.
       Thus,  the class of the Log-GARCH models generates a richer class of volatilities than the EGARCH.

\subsection{Strict stationarity}
We now show that the introduction of a
time varying intercept in the log-volatility  of Model (\ref{logGARCHSTABLE}) does not modify the
strict stationarity conditions of the Log-GARCH model.
The study % of the extended Log-GARCH model (\ref{logGARCHSTABLE})
being very similar to that of the Log-GARCH model (\ref{logGARCH}) in FWZ, details are omitted.
Let $\omega_t=\omega+\sum_{i=1}^q\omega_{i-}1_{\{\epsilon_{t-i}<0\}}$.
Because coefficients equal to zero can always be added, it is not restrictive to assume $p>1$ and $q>1$.
Let the vectors
\begin{eqnarray*}
\bepsilon_{t,q}^+&=&(1_{\{\epsilon_{t}>0\}}\log \epsilon_t^2,\dots,1_{\{\epsilon_{t-q+1}>0\}}\log \epsilon_{t-q+1}^2)'\in \mathbb{R}^q,\\
\bepsilon_{t,q}^-&=&(1_{\{\epsilon_{t}<0\}}\log \epsilon_t^2,\dots,1_{\{\epsilon_{t-q+1}<0\}}\log \epsilon_{t-q+1}^2)'\in \mathbb{R}^q,\\
\bz_t&=&(\bepsilon_{t,q}^+,\bepsilon_{t,q}^-,\log \sigma_t^2,\dots,\log \sigma_{t-p+1}^2)'\in \mathbb{R}^{2q+p},\\
\bb_t&=&\left((\omega_t+\log\eta_t^2)1_{\{\eta_{t}>0\}},\bzero'_{q-1},(\omega_t+\log\eta_t^2)1_{\{\eta_{t}<0\}},\bzero'_{q-1},\omega_t,\bzero'_{p-1}\right)'\in \mathbb{R}^{2q+p},
\end{eqnarray*} and the matrix
%\begin{equation}
%\label{A_tmulti}
$$
\bC_t=\left(\begin{array}{ccc}
1_{\{\eta_{t}>0\}}\balpha'_+&1_{\{\eta_{t}>0\}}\balpha'_-&1_{\{\eta_{t}>0\}}\bbeta'  \\
\begin{array}{cc}\bI_{q-1} &\bzero_{q-1}\end{array}  &\bzero_{(q-1)\times q}&\bzero_{(q-1)\times p}\\
1_{\{\eta_{t}<0\}}\balpha'_+&1_{\{\eta_{t}<0\}}\balpha'_-&1_{\{\eta_{t}<0\}}\bbeta'  \\
\bzero_{(q-1)\times q}&\begin{array}{cc}\bI_{q-1} &\bzero_{q-1}\end{array}   &\bzero_{(q-1)\times p}\\
\balpha'_+&\balpha'_-&\bbeta'  \\
\bzero_{(p-1)\times q}   &\bzero_{(p-1)\times q}&\begin{array}{cc}\bI_{p-1} &\bzero_{p-1}\end{array}\\
\end{array}\right).
$$%\end{equation}
Model (\ref{logGARCH}) is rewritten in matrix form as
$$%\begin{equation}\label{vectorform}
\bz_t=\bC_t\bz_{t-1}+\bb_t.$$%\end{equation}
%We have implicitly assumed $p>1$ and $q>1$ to write $\bC_t$ and  $\bb_t$, but obvious changes of notation can be employed when $p\leq 1$ or $q\leq 1$.
Let $\gamma ({\bf C})$
be the top Lyapunov exponent of the
sequence ${\bf C}=\{\bC_t, t\in \mathbb{Z}\}$,
$$
\gamma ({\bf C})=\lim_{t\to \infty} \frac{1}{t} E \left(\log
\|\bC_t\bC_{t-1}\ldots \bC_{1}\|\right) = \inf_{t\geq 1}\;
\frac{1}{t}\; E (\log \|\bC_t\bC_{t-1}\ldots \bC_{1}\|).
$$
%The choice of the norm is obviously unimportant for the value of the
%top Lyapunov exponent. However, in the sequel, the matrix norm will
%be assumed to be multiplicative.
It can be noted that the sequence $(\bC_t, \bb_t)$ is only strictly stationary and ergodic (not iid) but this property suffices
to extend the proof of Theorem 2.1 in FWZ.
%We then have the following extension of FWZ (2013).
\begin{theo}\label{theostatio}
Assume that $E\log^+|\log\eta_0^2|<\infty$.
A sufficient condition for the
existence of a strictly stationary solution
 to the AS-Log-GARCH model (\ref{logGARCHSTABLE}) is $\gamma ({\bf C})<0$.
%where $\gamma ({\bf C_0})$
%is the top Lyapunov exponent of the
%sequence ${\bf C_0}=\{C_t, t\in \mathbb{Z}\}$ defined in (\ref{A_tmulti}).
When $\gamma ({\bf C})<0$, there exists only one  stationary solution, which is non anticipative  and  ergodic.
\end{theo}
It follows that the presence of the coefficients $\omega_{i-}$ does not modify the stationarity condition.

\section{QML estimation of the AS-Log-GARCH model}
\label{sec2}
We turn to the inference of the AS-Log-GARCH model.
%Let $\bomega_{-}=(\omega_{1-},\dots,\omega_{q-})'$, $\balpha_{+}=(\alpha_{1+},\dots,\alpha_{q+})'$,
%$\balpha_{-}=(\alpha_{1-},\dots,\alpha_{q-})'$, and
%$\bbeta=(\beta_{1},\dots,\beta_{p})'$.
Let $\epsilon_1,\dots,\epsilon_n$ be observations of the stationary solution of (\ref{logGARCHSTABLE}),
where $\btheta=(\omega,\bomega_-',\balpha_+',\balpha_-',\bbeta')'$ is equal to an unknown value $\btheta_0$ belonging to some parameter space $\Theta\subset \mathbb{R}^d$, with $d=3q+p+1$.
A QMLE of $\btheta_0$ is defined as any measurable solution
$\widehat{\btheta}_n$  of
\begin{equation}\label{qml}
  \widehat{\btheta}_n=\argmin_{\btheta\in\Theta}\widetilde{Q}_n(\btheta),
\end{equation}
with
$$
  \widetilde{Q}_n(\btheta) =n^{-1}\sum_{t=r_0+1}^n\widetilde{\ell}_t(\btheta),\qquad
 \widetilde{\ell}_t(\btheta)= \frac{\epsilon_t^2}{\widetilde{\sigma}_t^2(\btheta)} +\log \widetilde{\sigma}_t^2(\btheta),
$$
where $r_0$ is a fixed integer and $\log \widetilde{\sigma}_t^2(\btheta)$ is recursively defined by
$\log \widetilde{\sigma}_t^2(\btheta)=\omega +\sum_{i=1}^q\left( \alpha_{i+}\log\epsilon_{t-i}^2 1_{\{\epsilon_{t-i}> 0\}}+(\omega_{i-}+ \alpha_{i-}\log\epsilon_{t-i}^2)1_{\{\epsilon_{t-i}<0\}}\right)
+\sum_{j=1}^p\beta_j\log \widetilde{\sigma}_{t-j}^2(\btheta),$ for $t=1,2,\dots,n$,
using initial values for $\epsilon_0,\dots, \epsilon_{1-q},\widetilde{\sigma}_{0}^2(\btheta),\dots, ,\widetilde{\sigma}_{1-p}^2(\btheta)$. We assume that
these initial values are such that
there exists a real random variable $K$ independent of $n$ satisfying
\begin{equation}
\label{condvi}
\sup_{\btheta\in\Theta}\left|\log\sigma_t^2(\btheta)-\log\widetilde{\sigma}_t^2(\btheta)\right|<K,\quad\mbox{a.s. for }t=q-p+1,\dots, q,
\end{equation}
where $\sigma_t^2(\btheta)$ is defined by
\begin{eqnarray}
\nonumber
\label{def:sig}
{\cal B}_{\btheta}(B)\log\sigma_t^2(\btheta)&=&\omega+{\cal O}^-_{\btheta}(B)1_{\{\epsilon_t<0\}}+
{\cal A}^+_{\btheta}(B)1_{\{\epsilon_{t}>0\}}\log\epsilon_{t}^2\\&&+
{\cal A}^-_{\btheta}(B)1_{\{\epsilon_{t}<0\}}\log\epsilon_{t}^2,
\end{eqnarray}
where   $B$ is the the lag operator and, for any $\btheta\in\Theta$, ${\cal A}^+_{\btheta}(z)=\sum_{i=1}^q\alpha_{i,+}z^{i}$,
${\cal A}^-_{\btheta}(z)=\sum_{i=1}^q\alpha_{i,-}z^{i}$, and
${\cal B}_{\btheta}(z)=1-\sum_{j=1}^p\beta_{j}z^{j}$ and ${\cal O}^-_{\btheta}(z)=\sum_{i=1}^q\omega_{i-}z^i$.
By
convention, ${\cal A}^{+}_{\btheta}(z)=0$, ${\cal A}^{-}_{\btheta}(z)=0$ and ${\cal O}^{-}_{\btheta}(z)=0$ if $q=0$, and  ${\cal
B}_{\btheta}(z)=1$ if $p=0$.
%To extend the results of FWZ (2013) one has to check  the identifiability of the extended version, see Remark \ref{ident?}. More precisely, we have to check the following identifiability property for the model \eqref{logGARCHSTABLE}:
%$$\sigma_1^2(\btheta)=\sigma_1^2(\btheta_0) \mbox{ a.s. }\mbox{ then }
%\btheta=\btheta_0.$$
%We also write ${\bf C}(\btheta_0)$ instead of ${\bf C}$ to emphasize  that the unknown parameter is $\btheta_0$.
Theorem~\ref{theostatio} shows that a strict stationarity condition of the Log-GARCH can be obtained from the behaviour of the sequence ${\bf C}$. As in FWZ, it can be shown that moment conditions can be obtained by constraining the matrix
\begin{align}
\label{MatA}
\bA_t&=
\left(\begin{array}{cccc}
\mu_1(\eta_{t-1})&\dots& \mu_{r-1}(\eta_{t-r+1})& \mu_r(\eta_{t-r})\\
\multicolumn{3}{c}{\bI_{r-1}} & \bzero_{r-1}
\end{array}\right),&
\end{align}
where $r=\max(p,q)$ and
$\mu_i(\eta_{t})=\alpha_{i+}1_{\{\eta_{t}>0\}}+\alpha_{i-}1_{\{\eta_{t}<0\}}+\beta_i$
with the convention $\alpha_{i+}=\alpha_{i-}=0$ for $i>p$ and
$\beta_{i}=0$ for $i>q$. The spectral radius of a square matrix $\bA$ is denoted by $\rho(\bA)$. For any vector or matrix $\bA$, we denote by
$\mathrm{Abs}({\bA})$ the matrix whose
elements are the absolute values of the corresponding elements of
$\bA$.

The following assumptions will be used to establish the strong consistency and asymptotic normality of the QMLE.
\begin{itemize}
\item[ {\bf A1:}]
\hspace*{1em} $ \btheta_0\in\Theta$  and  $\Theta$ is compact.
\item[ {\bf A2:}]
\hspace*{1em} $\gamma \left\{{\bf C}\right\}<0 \quad
$ and $\quad \forall \btheta\in\Theta,\quad |{\cal
B}_{\btheta}(z)|=0 \Rightarrow |z|>1.$
\item[ {\bf A3:}]
\hspace*{1em} %$P(\eta_0^2=0)=0$,
the support of $\eta_0$ contains at least two positive values and two negative values,
$E\eta_0^2=1$ and $E|\log\eta_0^2|^{s_0}<\infty$ for some
$s_0>0$.
\item[ {\bf A4:}] \hspace*{1em} If $p>0$ and $q>1$, there is no common root to the polynomials
${\cal O}^-_{\btheta_0}(z)$, ${\cal A}^+_{\btheta_0}(z)$, ${\cal A}^-_{\btheta_0}(z)$  and  ${\cal B}_{\btheta_0}(z)$. Moreover $(\bomega_{0-},\balpha_{0+},\balpha_{0-})\neq
0$ and  $|\omega_{0q-}||\alpha_{0q+}||\alpha_{0q-}|+|\beta_{0p}|\neq 0$ if $p>0$.
\item[ {\bf A5:}]
\hspace*{1em} $E\left|\log \epsilon_t^2\right|<\infty$.
\item[ {\bf A6:}]
\hspace*{1em} $\btheta_0\in\stackrel{\circ}{\Theta}$ and
$\kappa_4:=E(\eta_0^4)<\infty$.
\item[ {\bf A7:}]\hspace*{1em} There exists some $s_0>0$ such that $E\exp(s_0|\log\eta_0^2|)<\infty$ and $\rho\left\{\mathrm{ess}\sup \mathrm{Abs}({\bA}_{1})\right\}<1,$ where ${\bA}_{1}$ is defined by \eqref{MatA}.
\hspace*{1em}
\end{itemize}
In the case $p=q=1$, omitting the index $i$,  Assumption {\bf A2} simplifies to
the conditions
$|\alpha_{0+}+\beta_0|^a|\alpha_{0-}+\beta_0|^{1-a}<1$, where $a=P(\eta_0>0)$,
and $|\beta|<1, \forall \btheta\in\Theta$ (see FWZ, Example 2.1).

%Note that {\bf A3} is slightly stronger than the corresponding condition in FWZ (2013) that required only two positive and negative values.
Let $\nabla Q=(\nabla_1 Q,\dots, \nabla_d
Q)'$ and $\mathbb H Q=(\mathbb H_{1.} Q',\dots, \mathbb H_{d.}Q')'$
be the vector and matrix of the first-order and second-order partial
derivatives of a function $Q:\Theta\to \mathbb{R}$.
%\begin{theo}[{\bf Strong consistency of the QMLE}]\label{consistency}
%Let $(\widehat{\btheta}_n)$ be a sequence of QMLE satisfying (\ref{qml}), where the $\epsilon_t$'s follow the asymmetric Log-GARCH model of parameter $\btheta_0$.
%Under the assumptions (\ref{condvi}) and {\bf A1}-{\bf A5}, almost surely $\widehat{\btheta}_n\to \btheta_0$ as $n\to\infty$.
%\end{theo}
%\noindent{\bf Proof of Theorem \ref{consistency}.}
%The proof is similar than the one of Theorem 4.1 in FWZ (2013).
\begin{theo}[{\bf Asymptotic properties of the QMLE}]\label{normality}
Let $(\widehat{\btheta}_n)$ be a sequence of QMLE satisfying (\ref{qml}), where $(\epsilon_t)$ is the stationary solution of the AS-Log-GARCH model (\ref{logGARCHSTABLE}) with parameter $\btheta_0$.
Under the assumptions (\ref{condvi}) and {\bf A1}-{\bf A5}, $\widehat{\btheta}_n\to \btheta_0$ a.s. as $n\to\infty$.
If, moreover,
 {\bf A6}-{\bf A7} hold we have $\sqrt n (\widehat\btheta_n-\btheta_0)\stackrel{d}{\to}\mathcal N(\bzero,(\kappa_4-1)\bf J^{-1})$ as $n\to\infty$, where
 ${\bf J}=E[\nabla \log \sigma_t^2(\btheta_0) \nabla \log
\sigma_t^2(\btheta_0)']$ is a positive definite matrix and
 $\stackrel{d}{\to}$ denotes convergence in distribution.
\end{theo}
%\noindent{\bf Proof of Theorem \ref{normality}.}
\noindent{\bf Proof:}
The proof is similar to those of Theorems 4.1-4.2 of FWZ.
We will only
show the identifiability of the extended model, that is,
$$\sigma_1^2(\btheta)=\sigma_1^2(\btheta_0) \mbox{ a.s. } \quad \Rightarrow\quad
\btheta=\btheta_0.$$
Note that if the left-hand side holds, %$\log\sigma_1^2(\btheta)=\log\sigma_1^2(\btheta_0)$ a.s.,
by stationarity
we have $\log\sigma_t^2(\btheta)=\log\sigma_t^2(\btheta_0)$ for all
$t$. From the equality \eqref{def:sig} we then have, almost surely,
\begin{eqnarray*}&&\left\{\frac{{\cal O}^-_{\btheta}(B)}{{\cal
B}_{\btheta}(B)}-\frac{{\cal O}^-_{\btheta_0}(B)}{{\cal
B}_{\btheta_0}(B)}\right\}1_{\{\epsilon_{t}<0\}}+\left\{\frac{{\cal A}^+_{\btheta}(B)}{{\cal
B}_{\btheta}(B)}-\frac{{\cal A}^+_{\btheta_0}(B)}{{\cal
B}_{\btheta_0}(B)}\right\}1_{\{\epsilon_{t}>0\}}\log\epsilon_{t}^2\\
&&+
\left\{\frac{{\cal A}^-_{\btheta}(B)}{{\cal
B}_{\btheta}(B)}-\frac{{\cal A}^-_{\btheta_0}(B)}{{\cal
B}_{\btheta_0}(B)}\right\}1_{\{\epsilon_{t}<0\}}\log\epsilon_{t}^2=\frac{\omega_0}{{\cal
B}_{\btheta_0}(1)}-\frac{\omega}{{\cal B}_{\btheta}(1)}.\end{eqnarray*} Throughout the paper let
$R_t$ denote any generic random variable, whose value can be modified from one line to the other, which is measurable with respect to
$\sigma\left(\{\eta_u,u\leq t\}\right)$. If
\begin{equation}
\label{hypopasvraie}
\frac{{\cal O}^-_{\btheta}(B)}{{\cal B}_{\btheta}(B)}\neq \frac{{\cal O}^-_{\btheta_0}(B)}{{\cal B}_{\btheta_0}(B)}\;\mbox{ or }\;
\frac{{\cal A}^+_{\btheta}(B)}{{\cal B}_{\btheta}(B)}\neq \frac{{\cal A}^+_{\btheta_0}(B)}{{\cal B}_{\btheta_0}(B)}\;\mbox{ or }\;
\frac{{\cal A}^-_{\btheta}(B)}{{\cal B}_{\btheta}(B)}\neq \frac{{\cal A}^-_{\btheta_0}(B)}{{\cal B}_{\btheta_0}(B)},
\end{equation}
%\begin{equation}
%\label{hypopasvraie}
%\frac{{\cal O}^-_{\btheta}(B)}{{\cal B}_{\btheta}(B)}\neq \frac{{\cal O}^-_{\btheta_0}(B)}{{\cal B}_{\btheta_0}(B)}\quad\mbox{ or }\quad
%\frac{{\cal A}^+_{\btheta}(B)}{{\cal B}_{\btheta}(B)}\neq \frac{{\cal A}^+_{\btheta_0}(B)}{{\cal B}_{\btheta_0}(B)}\quad\mbox{ or }\quad
%\frac{{\cal A}^-_{\btheta}(B)}{{\cal B}_{\btheta}(B)}\neq \frac{{\cal A}^-_{\btheta_0}(B)}{{\cal B}_{\btheta_0}(B)},
%\end{equation}
 there exists a non null $(c_+,c_-,d_-)\in \mathbb{R}^3$,
 such that $$d_- 1_{\eta_t<0}+c_+ 1_{\{\eta_t> 0\}}\log\epsilon_t^2+c_- 1_{\{\eta_t<0\}}\log\epsilon_t^2+R_{t-1}=0\quad \mbox{ a.s.}$$  This is equivalent to the two equations
$$\left(c_+ \log\eta_t^2+c_+ \log\sigma_t^2+R_{t-1}\right)1_{\{\eta_t>0\}}=0\quad \mbox{ a.s.}$$  and $$\left(d_- +c_- \log\eta_t^2+c_- \log\sigma_t^2+R_{t-1}\right)1_{\{\eta_t<0\}}=0\quad \mbox{ a.s.}$$ Note that if an equation of the form $a\log x^2 1_{\{x>0\}}+b1_{\{x>0\}}=0$ admits two positive solutions then $a=0$. This result, {\bf A3}, and the independence between $\eta_t$ and $(\sigma_t^2,R_{t-1})$ imply that $c_+=0$ and $R_{t-1}=0$. Similarly we obtain $c_-=0$. Plugging $c_+=c_-=0$ in the equations above yields %$d_-=R_{t-1}=0$, which leads to
$c_+=c_-=d_-=0$ that is a contradiction. We conclude that (\ref{hypopasvraie}) cannot hold true, and the conclusion follows from {\bf A4}.
\zak

\section{Test of  AS-Log-GARCH}\label{sec4}
In this section, we are interested in testing the AS-Log-GARCH specification against more general formulations, including
both the Log-GARCH and the EGARCH models. For our testing problem, we therefore introduce the general model
\begin{equation}\label{logGARCHandEGARCH}
\left\{\begin{array}{lll}
\epsilon_t&=&\sigma_t\eta_t,\quad \\
\log \sigma_t^2&=&\omega_0+\sum_{i=1}^q\omega_{0,i-}1_{\{\epsilon_{t-i}<0\}}\\&&+\sum_{i=1}^q\left(\alpha_{0,i+}1_{\{\epsilon_{t-i}>0\}}+\alpha_{0,i-}1_{\{\epsilon_{t-i}<0\}}\right)\log\epsilon_{t-i}^2\\
&&+\sum_{j=1}^p\beta_{0j}\log \sigma_{t-j}^2
+ \sum_{k=1}^{\ell} \gamma_{0,k+} \eta_{t-k}^++ \gamma_{0,k-} \eta_{t-k}^-.
\end{array}\right.
\end{equation}
%where $x^+=\max (x,0), x^-=\max (-x,0)$, and  under the same
%assumptions on the sequence $(\eta_t)$ as in Model (\ref{logGARCH}).
Let $\bvartheta_0=(\btheta_0', \bgamma_0')'$ where
$\bgamma_0=(\gamma_{01,+}, \gamma_{01,-},\dots,\gamma_{0\ell,-})'$
and $\btheta_0$ is as in Section~\ref{sec2}.

We wish to test the hypothesis that, in (\ref{logGARCHandEGARCH}),
$$H^{\bgamma}_0: \bgamma_0=\bzero_{2\ell\times 1}\quad\mbox{against}\quad H_1^{\bgamma}:
\bgamma_0\neq \bzero_{2\ell\times 1}.$$ In the time series literature, similar testing
problems are solved by a standard test, using for example the Wald,
Lagrange-Mutiplier (LM) or Likelihood-Ratio (LR) principle. See
among others Luukkonen, Saikkonen and Teräsvirta (1988), Francq,
Horv\'ath and Zakoïan (2010).

A difficulty, in the present framework, is that we do not have a
consistent estimator of the parameter $\bvartheta_0$. Two problems
arise to prove that the QMLE is consistent. First, the stationarity
conditions of Model (\ref{logGARCHandEGARCH}) are unknown. Second,
due to the presence of the $|\eta_{t-k}|$'s, it seems extremely
difficult to obtain invertibility conditions allowing to write $\log
\sigma_t^2 (\bvartheta)$ (where $\bvartheta$ denotes any parameter
value) as a function of the observations.

To circumvent these problems, we propose a LM approach. Denote by
$\widehat{\bvartheta}_n^c$ the constrained (by $H_0^{\bgamma}$)
estimator of $\bvartheta_0$, defined by
$$%\begin{equation}\label{qmlconst}
  \widehat{\bvartheta}_n^c=  (\widehat{\btheta}_n', \mathbf{0}_{1\times 2\ell})'
%  \argmin_{\bvartheta\in\Theta}\widetilde{Q}_n(\btheta),
$$%\end{equation}
where $\widehat{\btheta}_n$ is the QMLE of the AS-Log-GARCH parameters
defined in (\ref{qml}).

For any $\bvartheta$ in $\Theta \times \mathbb{R}^{2\ell}$, define
$\log \widetilde{\sigma}_t^2(\bvartheta)$ recursively, for
$t=1,2,\dots,n$, by \begin{eqnarray*}\log
\widetilde{\sigma}_t^2(\bvartheta)&=&\omega+\sum_{i=1}^q\omega_{i-}1_{\{\epsilon_{t-i}<0\}}+\sum_{i=1}^q\left(\alpha_{i+}1_{\{\epsilon_{t-i}>0\}}+\alpha_{i-}1_{\{\epsilon_{t-i}<0\}}\right)\log\epsilon_{t-i}^2
\\&&+\sum_{j=1}^p\beta_j\log \widetilde{\sigma}_{t-j}^2(\bvartheta)
+ \sum_{k=1}^{\ell} (\gamma_{k+} \epsilon_{t-k}^++\gamma_{k-}
\epsilon_{t-k}^-) e^{-\frac12\log
\widetilde{\sigma}_{t-k}^2(\bvartheta)},
\end{eqnarray*}
using positive
initial values for $\epsilon_0^2,\dots,
\epsilon_{1-\max
(q,\ell)}^2,\widetilde{\sigma}_{0}^2(\bvartheta),\dots,
,\widetilde{\sigma}_{1-\max (p,\ell)}^2(\bvartheta)$. The  random
vector  $\frac{\partial}{\partial \bvartheta}
\log\widetilde{\sigma}_t^2(\bvartheta)$ satisfies
\begin{eqnarray*}\frac{\partial}{\partial \bvartheta} \log\widetilde{\sigma}_t^2(\bvartheta)&=&
-\frac12 \sum_{k=1}^{\ell} (\gamma_{k+} \epsilon_{t-k}^++\gamma_{k-}
\epsilon_{t-k}^-) e^{-\frac12\log
\widetilde{\sigma}_{t-k}^2(\bvartheta)}\frac{\partial}{\partial \bvartheta} \log\widetilde{\sigma}_{t-k}^2(\bvartheta)\\
&&+\sum_{j=1}^p\beta_j \frac{\partial}{\partial \bvartheta} \log\widetilde{\sigma}_{t-j}^2(\bvartheta)+
\left(\begin{array}{c}1\\{\bf 1}_{t-1,q}^-\\
\bepsilon_{t-1,q}^+\\\bepsilon_{t-1,q}^-\\
\widetilde{\bsigma}^2_{t-1,p}(\bvartheta)\\ \widetilde{\boldeta}_{t-1}(\bvartheta)
\end{array}\right),
\end{eqnarray*}
where
\begin{eqnarray*}
\widetilde{\bsigma}^2_{t,p}(\bvartheta)&=&(\log
\widetilde{\sigma}_{t}^2(\bvartheta), \ldots, \log
\widetilde{\sigma}_{t-p+1}^2(\bvartheta))',\\ \widetilde{\boldeta}_{t}(\bvartheta)&=&(\epsilon_{t}^+ e^{-\frac12\log
\widetilde{\sigma}_{t}^2(\bvartheta)}, \epsilon_{t}^- e^{-\frac12\log
\widetilde{\sigma}_{t}^2(\bvartheta)}, \ldots, \epsilon_{t-\ell+1}^-
e^{-\frac12\log \widetilde{\sigma}_{t-\ell+1}^2(\bvartheta)})'.
\end{eqnarray*}
With a slight abuse of notation we write
$\widetilde{\sigma}_t^2(\bvartheta)=\widetilde{\sigma}_t^2(\btheta)$
when $\bvartheta=(\btheta', \mathbf{0}_{1\times 2\ell})'$, that is
when $\bvartheta$ satisfies $H_0^{\bgamma}$. Similarly, to avoid
introducing new notations we still define the criterion function by
$$
  \widetilde{Q}_n(\bvartheta) =n^{-1}\sum_{t=r_0+1}^n\widetilde{\ell}_t(\bvartheta),\qquad \mbox{where} \quad
 \widetilde{\ell}_t(\bvartheta)= \frac{\epsilon_t^2}{\widetilde{\sigma}_t^2(\bvartheta)} +\log \widetilde{\sigma}_t^2(\bvartheta).
$$
%We a slight abuse of notation we also let ${\cal
%B}_{\bvartheta}={\cal B}_{\btheta}$.
To derive a LM test, we need to
find the asymptotic distribution of
\begin{eqnarray*}
\frac1{\sqrt{n}}\sum_{t=1}^n \frac{\partial}{\partial \bvartheta} \widetilde{\ell}_t(\widehat{\bvartheta}_n^c)=
\left(\begin{array}{c}
\mathbf{0}_{d\times 1}\\\mathbf{S}_n:=
\frac1{\sqrt{n}}\sum_{t=1}^n (1-\widehat{\eta}_t^2)\widehat{\bnu}_t
%{\cal B}_{\widehat{\bvartheta}_n^c}^{-1}(B)\widehat{\boldeta}_{t-1,\ell}
%\widetilde{\by}_{t-1,\ell}(\widehat{\bvartheta}_n^c)
\end{array}\right),\quad \widehat{\bnu}_t={\cal B}_{\widehat{\btheta}_n}^{-1}(B)\widehat{\boldeta}_{t-1}
\end{eqnarray*}
where
$\widehat{\eta}_t=\epsilon_{t}/\widetilde{\sigma}_t(\widehat{\btheta}_n)$ for $t\geq 1$,
$\widehat{\eta}_t=0$ for $t\leq 0$,
%with by convention $\widetilde{\sigma}_t(\widehat{\btheta}_n)$ for
%$t\leq 0$
and $\widehat{\boldeta}_{t}=(\widehat{\eta}_t^+,
\widehat{\eta}_t^-, \ldots, \widehat{\eta}_{t-\ell+1}^-)'.$ Note that the
nullity of the first $d$ components of the score follows from the definition of
$\widehat{\bvartheta}_n^c$ as a maximizer of the quasi-likelihood in
the restricted model. The invertibility of the lag polynomial ${\cal
B}_{\widehat{\btheta}_n}(B)$ follows from {\bf A2}.

The following quantities are used to define the LM test statistic.
Recall that $\nabla$ denotes the differentiation operator with
respect to the components of $\btheta$. Let
\begin{eqnarray*}
\widehat{\bcalJ}_{11} &=& \frac1n \sum_{t=1}^n
\widehat{\bnu}_t\widehat{\bnu}_t'-\left(\frac1n \sum_{t=1}^n
\widehat{\bnu}_t\right)\left(\frac1n \sum_{t=1}^n
\widehat{\bnu}_t'\right),
\quad \widehat{\kappa}_4 -1=\frac1n
\sum_{t=1}^n (1-\widehat{\eta}_t^2)^2,\\
\widehat{\bf J}&=&\frac1n \sum_{t=1}^n \nabla \log
\widetilde{\sigma}_t^2(\widehat{\btheta}_n)\nabla' \log
\widetilde{\sigma}_t^2(\widehat{\btheta}_n),\qquad
\widehat{\bOmega}=\frac1n \sum_{t=1}^n \widehat{\bnu}_t\nabla' \log
\widetilde{\sigma}_t^2(\widehat{\btheta}_n),
\\
\widehat{\bcalJ}_{12} &=& -\left\{\widehat{\bOmega}-\left(\frac1n \sum_{t=1}^n
\widehat{\bnu}_t\right)\left(\frac1n \sum_{t=1}^n
\nabla' \log
\widetilde{\sigma}_t^2(\widehat{\btheta}_n)\right)\right\}\widehat{\bf J}^{-1}=\widehat{\bcalJ}_{21}',
\end{eqnarray*}
 and
$$\widehat{\bcalI}=\widehat{\bcalJ}_{11}+ \widehat\bOmega \widehat{\bf J}^{-1}\widehat\bOmega'+ \widehat{\bcalJ}_{12}\widehat\bOmega'+ \widehat\bOmega \widehat{\bcalJ}_{21}.$$
To derive the test, we need to slightly reinforce {\bf A3}
concerning the support of the distribution of $\eta_t$.
\begin{itemize}
\item[ {\bf A8:}] The support of $\eta_0$ contains at least three positive values and three negative values.
\end{itemize}
\begin{theo}[{\bf Asymptotic distribution of the LM test under $H_0^{\bgamma}$}]\label{test}
Under the assumptions of Theorem \ref{normality} (thus under $H_0^{\bgamma}$) and {\bf A8},
the matrix $\widehat{\bcalI}$ converges in
probability to a positive definite matrix ${\bcalI}$ and we have
$$\mathbf{LM}_n^{\bgamma}=(\widehat{\kappa}_4 -1)^{-1}\mathbf{S}_n'\widehat{\bcalI}^{-1}\mathbf{S}_n\stackrel{d}{\to} \chi^2_{2\ell}$$
where $\chi^2_{2\ell}$ denotes the chi-square distribution with
$2\ell$ degrees of freedom.
\end{theo}
%$\mathrm{\bf{LM}}_n(\gamma)=n\frac{\widetilde{\sigma}^2-\widehat{\sigma}_{\gamma}^2}{\widetilde{\sigma}^2}$
Denoting by $\chi_{\ell}^2(\alpha)$ the $\alpha$-quantile of the chi-square distribution with $\ell$ degrees of freedom, the AS-Log-GARCH$(p,q)$ model (\ref{logGARCHSTABLE}) is then rejected at the asymptotic level $\alpha$
when
%\begin{equation}
%\label{portmanteautest}
$\left\{\mathbf{LM}_n^{\bgamma}>\chi_{2\ell}^2(1-\alpha)\right\}.$
%\end{equation}

\noindent {\bf Proof:} %of Theorem  \ref{test}
For any $\bvartheta^c= (\btheta',
\mathbf{0}_{1\times 2\ell})'\in \Theta \times \{0\}^{2\ell}$, let
%${\sigma}_t^2(\bvartheta^c)$ denote the strictly stationary solution
%of
%\begin{eqnarray*}\log
%{\sigma}_t^2(\bvartheta^c)&=&\omega+\sum_{i=1}^q\left(\alpha_{i+}1_{\{\epsilon_{t-i}>0\}}+\alpha_{i-}1_{\{\epsilon_{t-i}<0\}}\right)\log\epsilon_{t-i}^2
%+\sum_{j=1}^p\beta_j\log {\sigma}_{t-j}^2(\bvartheta^c)
%\end{eqnarray*}
%and let
$
\eta_t(\btheta)=\frac{\epsilon_t}{{\sigma}_t(\btheta)},$
$\boldeta_{t}(\btheta)=({\eta}_t^+(\btheta), {\eta}_t^-(\btheta),\ldots,
{\eta}_{t-\ell+1}^-(\btheta))',$ $\bnu_t(\btheta)={\cal
B}_{\btheta}^{-1}(B)\boldeta_{t-1}(\btheta)$
%\widetilde{\eta}_t(\bvartheta^c)=\frac{\epsilon_t}{\widetilde{\sigma}_t(\btheta)}.$$
and let  $\widetilde{\eta}_t(\btheta), \widetilde{\boldeta}_{t}(\btheta),
\widetilde{\bnu}_t(\btheta)$ denote the corresponding quantities when
${\sigma}_t(\btheta)$ is replaced by $\widetilde{\sigma}_t(\btheta)$.
 Let also
$$\mathbf{S}_n(\btheta)= \frac1{\sqrt{n}}\sum_{t=1}^n
\{1-{\eta}_t^2(\btheta)\}{\bnu}_t(\btheta), \quad
\widetilde{\mathbf{S}}_n(\btheta)= \frac1{\sqrt{n}}\sum_{t=1}^n
\{1-\widetilde{\eta}_t^2(\btheta)\}\widetilde{\bnu}_t(\btheta).$$ Let
$\mathbf{S}_{n,i}$ denote the $i$-th component of
$\mathbf{S}_n=\widetilde{\mathbf{S}}_n(\widehat{\btheta}_n)$, for $i=1,
\ldots 2\ell$.
 A Taylor
expansion gives, for some ${\btheta}_*$ between
$\widehat{\btheta}_n$ and ${\btheta}_0$,
\begin{eqnarray}\label{4.3}
\mathbf{S}_{n,i} = %=\widetilde{\mathbf{S}}_n(\widehat{\btheta}_n)=
\widetilde{\mathbf{S}}_{n,i} ({\btheta}_0)+\frac1{\sqrt{n}}\frac{\partial \widetilde{\mathbf{S}}_{n,i}}{\partial \btheta'}({\btheta}_*)
\sqrt{n}(\widehat{\btheta}_n-
{\btheta}_0).
\end{eqnarray}
Recall that ${\bf J}=E[\nabla \log \sigma_t^2(\btheta_0) \nabla \log
\sigma_t^2(\btheta_0)']$ and define \begin{eqnarray*} \boldsymbol{\bcalJ}&=&\left(\begin{array}{cc}{\bcalJ}_{11}&{\bcalJ}_{12}\\{\bcalJ}_{21}&{\bcalJ}_{22}\end{array}\right),\quad \mbox{where } \quad {\bcalJ}_{11}=\mbox{Var}\{\bnu_t(\btheta_0)\},\quad {\bcalJ}_{22}={\bf
J}^{-1}\\ &&\;\;\;\qquad \qquad {\bcalJ}_{12}={\bcalJ}_{21}'=-\mbox{Cov}\{\bnu_t(\btheta_0), \nabla \log
\sigma_t^2(\btheta_0)\} {\bf J}^{-1},  \end{eqnarray*} and
$${\bcalI}={\bcalJ}_{11}+ \bOmega{\bf J}^{-1}\bOmega'+ {\bcalJ}_{12}\bOmega'+ \bOmega{\bcalJ}_{21},$$
where $\bOmega=E\{\bnu_t(\btheta_0)\nabla' \log
\sigma_t^2(\btheta_0)\}.$ Let $\bOmega_i=E\{\nu_{t,i}(\btheta_0)\nabla' \log
\sigma_t^2(\btheta_0)\}$, where
$\nu_{t,i}(\btheta_0)$ denotes the $i$-th component of
$\bnu_{t}(\btheta_0)$, for $i=1, \ldots 2\ell$.

The advanced result
 is obtained by showing the following intermediate steps: under
 $H_0^{\bgamma}$, as $n\to \infty$,
\begin{eqnarray*}
&&i) \;  \sup_{\btheta \in {\cal
V}(\btheta_0)}\left\|{\mathbf{S}}_n(\btheta)-
\widetilde{\mathbf{S}}_n(\btheta)
\right\|\to
0, \quad \sup_{\btheta \in {\cal
V}(\btheta_0)}\frac1{\sqrt{n}}\left\|\frac{\partial {\mathbf{S}}_n}{\partial \btheta'}(\btheta)-
\frac{\partial \widetilde{\mathbf{S}}_n}{\partial \btheta'}(\btheta)
\right\|\to
0,\\&& \qquad
\;\mbox{ in probability, }\\
&& ii) \; \left(\begin{array}{c}{\mathbf{S}}_n({\btheta}_0) \\ \sqrt{n}(\widehat{\btheta}_n-
{\btheta}_0)\end{array}\right)\stackrel{d}{\to} {\cal
N} (\mathbf{0},(\kappa_4 -1){\bcalJ}), %\quad \mbox{where}
\\&& iii)\;\mbox{ There exists a neighborhood of ${\cal V}(\btheta_0)$ of $\btheta_0$, such that, for } i=1, \ldots , \ell\\
&&\qquad E
\sup_{\btheta \in {\cal
V}(\btheta_0)}\left\|\mathbb H
\left[\{1-\eta_t^2(\btheta)\}{\cal
B}_{\btheta}^{-1}(B)|\eta_{t-i}(\btheta)|\right]\right\|<\infty,
 \\&& iv)
\;  \frac1{\sqrt{n}} \frac{\partial {\mathbf{S}}_{n,i}}{\partial \btheta'}(\btheta_*)\to
\bOmega_i, %=E\{\nabla' \log \sigma_t^2(\btheta_0)\bnu_t(\bvartheta_0^c)\}
\qquad
\;\mbox{ in probability as $n\to \infty$,}
 \\
&& v) \; {\bcalI} \;\mbox{ is non-singular.}
%&& vi) \; \widehat{\bcalI} \to {\bcalI} \;\mbox{ in probability.}
\end{eqnarray*}
We will use the following Lemma, whose proof is similar to that of Lemma 4.2 in FWZ and is thus omitted.
\begin{lem}\label{lem:mom}
Under the assumptions of Theorem
\ref{normality}, for any $m>0$ there exists a neighborhood $\mathcal V $ of $\btheta_0$ such that $E[\sup_{\mathcal V } (\sigma_t^2/\sigma_t^2(\btheta))^m]<\infty$ and $E[\sup_{\mathcal V } |\log\sigma_t^2(\btheta)|^m]<\infty$.
\end{lem}
To prove the first convergence in {\it i)}, note that
\begin{eqnarray*}
\left\|{\mathbf{S}}_n(\btheta)- \widetilde{\mathbf{S}}_n(\btheta)
\right\|&\leq &\frac1{\sqrt{n}}
\sum_{t=1}^n |1-{\eta}_t^2(\btheta)|\left\|{\bnu}_t(\btheta)-\widetilde{\bnu}_t(\btheta)\right\|
\\
&&+\frac1{\sqrt{n}}
\sum_{t=1}^n |\widetilde{\eta}_t^2(\btheta)-{\eta}_t^2(\btheta)|\left\|\widetilde{\bnu}_t(\btheta)\right\| =S_1(\btheta)+S_2(\btheta).
\end{eqnarray*}
%As in FZ, using (\ref{condvi}) one can
We will show that there exist  $K>0$ and $ \rho\in (0,1)$, such that for almost all trajectories and  for all $\btheta\in\Theta,$
\begin{equation}\label{eq1}
\left|\frac{1}{\sigma_t^2(\btheta)}-\frac{1}{\widetilde{\sigma}_t^2(\btheta)}\right|\leq \frac{K\rho^t}{\sigma_t^2(\btheta)}.
\end{equation}
Similarly to the proof of (7.8) in FWZ, it can be shown that
$$
 \sup_{\btheta\in\Theta}\frac1{t}\log\left|\frac{1}{\sigma_t^2(\btheta)}-\frac{1}{\widetilde{\sigma}_t^2(\btheta)}\right|\leq \frac{a_{1t}}t+a_{2t},
$$
where $E|a_{1t}|<\infty$ and $\limsup_{t\to \infty} a_{2t}=\log \tilde{\rho}$
for some $\tilde{\rho}\in (0,1)$. We thus have
$$
\frac1{t}\log \sigma_t^2(\btheta)\left|\frac{1}{\sigma_t^2(\btheta)}-\frac{1}{\widetilde{\sigma}_t^2(\btheta)}\right|\leq \frac{\log \sigma_t^2(\btheta)}{t}+\frac{a_{1t}}t+a_{2t}.
$$
%using the fact that $1/{t}\log \sigma_t^2(\btheta)\to 0$ a.s. as $t\to \infty$. The latter convergence
The first term in the right-hand side converges a.s. to zero
as a consequence of Lemma 7.2 in FWZ and $E\sup_{\btheta\in\Theta}|\log \sigma_t^2(\btheta)|<\infty$, which follows from {\bf A5}. Thus (\ref{eq1}) is established.
Then we obtain
\begin{eqnarray*}
|\widetilde{\eta}_t^2(\btheta)-{\eta}_t^2(\btheta)|=\epsilon_t^2
\left|\frac1{\widetilde{\sigma}_t^2(\btheta)}-\frac1{{\sigma}_t^2(\btheta)}\right|\leq \epsilon_t^2\frac{K\rho^t}{\sigma_t^2(\btheta)}.
\end{eqnarray*}
Lemma~\ref{lem:mom} and the $c_r$ and H\"older inequalities  entail
that for sufficiently small $s\in(0,1)$, there exists a neighborhood
$\mathcal V $ of $\btheta_0$ such that
\begin{eqnarray*} E\sup_{\btheta \in {\cal
V}(\btheta_0)}S_2^s(\btheta)&\leq&
\frac{K}{n^{s/2}} \sum_{t=1}^n
\rho^{st}E\left[|\eta_t|^{2s}\sup_{\btheta \in {\cal
V}(\btheta_0)}\left\{\frac{\sigma_t^{2s}}{\sigma_t^{2s}(\btheta)}
\left\|\widetilde{\bnu}_t(\btheta)\right\|^s\right\}\right]\\
&\leq&
\frac{K}{n^{s/2}} \sum_{t=1}^n
\rho^{st}\to 0
\end{eqnarray*}
as $n\to \infty$. This entails $\sup_{\btheta \in {\cal
V}(\btheta_0)}S_2(\btheta)=o_P(1).$ Similarly, we have
$\sup_{\btheta \in {\cal V}(\btheta_0)}S_1(\btheta)=o_P(1).$ The
first convergence in {\it i)} follows and the second one is obtained
by the same arguments.

 To prove {\it ii)}, note that
\begin{eqnarray*}
%\left(\begin{array}{c}{\mathbf{S}}_n({\bvartheta}_0^c) \\
%\sqrt{n}(\widehat{\bvartheta}_n^c-
%{\bvartheta}_0^c)\end{array}\right)=
\left(\begin{array}{c}{\mathbf{S}}_n({\btheta}_0) \\
\sqrt{n}(\widehat{\btheta}_n-
{\btheta}_0)\end{array}\right)=
\frac1{\sqrt{n}}\sum_{t=1}^n
(1-{\eta}_t^2)\left(\begin{array}{c}\bnu_t(\btheta_0)\\
-{\bf J}^{-1}\nabla \log \sigma_t^2(\btheta_0)
\end{array}\right)+o_P(1).
\end{eqnarray*}
The convergence in distribution thus follows from the central limit
theorem for martingale differences.
% To show that ${\cal
%J}$ is non-singular, let $x=(x_i)\in \mathbb{R}^l$ such that
%$x'{\bcalJ}x=0$. It follows that $x'\bnu_t=0$ a.s. If $x_1=0$, we
%thus have
%$$x_1{\cal B}_{\bvartheta_0}^{-1}(B)|\eta_{t-1}|=\sum_{i=2}^{\ell}
%x_i{\cal B}_{\bvartheta_0}^{-1}(B)|\eta_{t-i}|, \quad a.s.$$ from
%which we deduce that $|\eta_{t-1}|$ is a linear combination of the
%$|\eta_{t-j}|$'s for $j>1$. This is in contradiction with
%$\mbox{Var}|\eta_t|>0$. Proceeding by induction we thus show that
%$x=0$, proving that
% ${\cal
%J}$ is non-singular.

To prove {\it iii)}, write ${\cal
B}_{\btheta}^{-1}(B)=\sum_{j=0}^{\infty}c_j(\btheta)B^j$.
We have
\begin{eqnarray*}
\mathbb H\{\eta_t(\btheta)\}&=&\eta_t(\btheta)\left[\frac14 \nabla \log\sigma_t^2(\btheta)\nabla' \log\sigma_t^2(\btheta)-
\frac12\mathbb H\{\log\sigma_t^2(\btheta)\}\right],\\
\mathbb H\{1-\eta_t^2(\btheta)\}&=&\eta_t^2(\btheta)\left[-\nabla \log\sigma_t^2(\btheta)\nabla' \log\sigma_t^2(\btheta)+
\mathbb H\{\log\sigma_t^2(\btheta)\}\right].
\end{eqnarray*}
It follows that, dropping  temporarily the term "$(\theta)$" to lighten the notation,
%\begin{eqnarray*}
%&&\mathbb H
%\{1-\eta_t^2(\btheta)\}c_j(\btheta)|\eta_{t-i-j}(\btheta)|\\&=&
%\eta_t^2(\btheta)\{-\nabla \log\sigma_t^2(\btheta)\nabla' \log\sigma_t^2(\btheta) + \mathbb H \log\sigma_t^2(\btheta)\}c_j(\btheta)|\eta_{t-i-j}(\btheta)| \\
%&&+
%\{1-\eta_t^2(\btheta)\}\{\mathbb Hc_j(\btheta)\}|\eta_{t-i-j}(\btheta)|\\
%&&+
%\{1-\eta_t^2(\btheta)\}c_j(\btheta)|\eta_{t-i-j}(\btheta)|\{\frac14 \nabla \log\sigma_t^2(\btheta)\nabla' \log\sigma_t^2(\btheta)
%-\frac12 \mathbb H \log\sigma_t^2(\btheta)\}\\
%&&
%+\eta_t^2(\btheta)\{\nabla \log\sigma_t^2(\btheta)\nabla'c_j(\btheta)+\nabla c_j(\btheta)\nabla'\log\sigma_t^2(\btheta)\}|\eta_{t-i-j}(\btheta)|\\
%&&-\frac12 \eta_t^2(\btheta)\{\nabla \log\sigma_t^2(\btheta)\nabla'\log\sigma_{t-i-j}^2(\btheta)
%+\nabla \log\sigma_{t-i-j}^2(\btheta)\nabla'\log\sigma_t^2(\btheta)\}|\eta_{t-i-j}(\btheta)|c_j(\btheta)\\
%&&-\frac12 \{1-\eta_t^2(\btheta)\}\{\nabla c_j(\btheta)\nabla'\log\sigma_{t-i-j}^2(\btheta)+\nabla \log\sigma_{t-i-j}^2(\btheta)\nabla'c_j(\btheta)\}|\eta_{t-i-j}(\btheta)|.
%\end{eqnarray*}
\begin{eqnarray*}
&&\mathbb H
\{1-\eta_t^2\}c_j|\eta_{t-i-j}|\\&=&
\eta_t^2\{-\nabla \log\sigma_t^2\nabla' \log\sigma_t^2 + \mathbb H \log\sigma_t^2\}c_j|\eta_{t-i-j}| \\
&&+
\{1-\eta_t^2\}\{\mathbb Hc_j\}|\eta_{t-i-j}|\\
&&+
\{1-\eta_t^2\}c_j|\eta_{t-i-j}|\{\frac14 \nabla \log\sigma_t^2\nabla' \log\sigma_t^2
-\frac12 \mathbb H \log\sigma_t^2\}\\
&&
+\eta_t^2\{\nabla \log\sigma_t^2\nabla'c_j+\nabla c_j\nabla'\log\sigma_t^2\}|\eta_{t-i-j}|\\
&&-\frac12 \eta_t^2\{\nabla \log\sigma_t^2\nabla'\log\sigma_{t-i-j}^2
+\nabla \log\sigma_{t-i-j}^2\nabla'\log\sigma_t^2\}|\eta_{t-i-j}|c_j\\
&&-\frac12 \{1-\eta_t^2\}\{\nabla c_j\nabla'\log\sigma_{t-i-j}^2+\nabla \log\sigma_{t-i-j}^2\nabla'c_j\}|\eta_{t-i-j}|.
\end{eqnarray*}
In view of Lemma \ref{lem:mom}, since
$\eta_t(\btheta)=\eta_t\sigma_t(\btheta_0)/\sigma_t(\btheta)$,
because $\nabla \log\sigma_t^2(\btheta)$   admits moments of any
order, and using the H\"older inequality, the conclusion follows.

To prove {\it iv)}, consider the following Taylor expansion about
$\btheta_0$
$$\frac1{\sqrt{n}} \frac{\partial {\mathbf{S}}_{n,i}}{\partial \btheta}(\btheta_*)=
\frac1{\sqrt{n}} \frac{\partial {\mathbf{S}}_{n,i}}{\partial
\btheta}(\btheta_0)+ \frac1{\sqrt{n}} \frac{\partial^2
{\mathbf{S}}_{n,i}}{\partial \btheta\partial
\btheta'}(\btheta^*)(\btheta_*-\btheta_0)$$ where $\btheta^*$ is
between $\btheta_*$ and $\btheta_0$. The a.s. convergence of
$\btheta^*$ to $\btheta_0$, {\it iii)} and the ergodic theorem imply
that, for $i=2k+1$  and for some neighborhood of $\btheta_0$
\begin{eqnarray*}
&&\lim\sup_{n\to \infty} \left\|\frac1{\sqrt{n}} \frac{\partial^2
{\mathbf{S}}_{n,i}}{\partial \btheta\partial
\btheta'}(\btheta^*)\right\|\\&\leq&
\lim\sup_{n\to \infty} \frac1n \sum_{t=1}^n
\sup_{\btheta \in {\cal
V}(\btheta_0)}\left\|\frac{\partial^2}{\partial \btheta\partial
\btheta'}
\{1-\eta_t^2(\btheta)\}{\cal
B}_{\btheta}^{-1}(B)\eta_{t-k-1}^+(\btheta)\right\|\\
&=&
E
\sup_{\btheta \in {\cal
V}(\btheta_0)}\left\|\frac{\partial^2}{\partial \btheta\partial
\btheta'}
\{1-\eta_t^2(\btheta)\}{\cal
B}_{\btheta}^{-1}(B)\eta_{t-k-1}^+(\btheta)\right\|<\infty.
\end{eqnarray*}
The same argument obviously applies for $i=2k$ and the conclusion
follows.

To prove {\it v)}, in view of (\ref{4.3}), it suffices to show that
${\bcalJ}$ is non-singular. %and $\bOmega$ has full row rank.
Suppose
there exist ${\bf x}=(x_i)\in \mathbb{R}^{2\ell}$ and ${\bf y}\in
\mathbb{R}^d$ such that
$$%\begin{equation}\label{eqident}
{\bf x}'\bnu_t(\btheta_0)+{\bf y}' {\bf J}^{-1}\nabla \log
\sigma_t^2(\btheta_0)=0, \quad \mbox{a.s.}
$$%\end{equation}
Recall that, in view of (\ref{def:sig}),
\begin{equation}
\nabla \log\sigma_t^2(\btheta_0)={\cal
B}_{\btheta_0}^{-1}(B)\left(1,{\bf 1}^{-'}_{t-1,q}, \bepsilon_{t-1,q}^{+'},
\bepsilon_{t-1,q}^{-'}, \bsigma^{2'}_{t-1,p}(\btheta_0)\right)'.\nonumber
\end{equation}
Letting ${\bf z}={\bf J}^{-1}{\bf y}=(z_i)$, we find that,
$x_1\eta_{t-1}^++
x_2\eta_{t-1}^-+z_21_{\{\eta_{t-1}>0\}}+z_{2+q}1_{\{\eta_{t-1}>0\}}\log \epsilon_{t-1}^2+
z_{2+2q}1_{\{\eta_{t-1}<0\}}\log \epsilon_{t-1}^2= R_{t-2}, \quad \mbox{a.s.}$
% where $R_{t-2}$ is a random variable belonging to the sigma-field
%$\sigma\left(\{\eta_u,u\leq t-2\}\right)$.
 Conditionally on $\eta_{t-1}>0$ we thus have
$$x_1\eta_{t-1}+z_2+
z_{2+q}\log \eta_{t-1}^2+z_{2+q}\log \sigma_{t-1}^2=
 R_{t-2}, \quad \mbox{a.s.}$$
% Similarly, conditional on $\eta_{t-1}<0$,
%$$x_2\eta_{t-1}
%=z_{2+q}\log \eta_{t-1}^2+z_{2+q}\log \sigma_{t-1}^2=
% R_{t-2}.$$
 By {\bf A8}, we find  $x_1=
 z_{2+q}=0.$ By conditioning on $\eta_{t-1}<0$, we similarly get $x_2=
 z_{2+2q}=0.$
 Thus
 $z_21_{\{\eta_{t-1}>0\}}= R_{t-2}, \quad \mbox{a.s.},$
 from which we deduce $z_2=R_{t-2}=0\quad \mbox{a.s.}$
%  We thus have, by (\ref{eqident}),
%conditionally on $\eta_{t-2}>0$,
%\begin{eqnarray*}
%x_3\eta_{t-2} &=&z_3\log \eta_{t-2}^2+z_{d+1}\log \sigma_{t-1}^2+
% h_{t-3}\\
%&=&(z_3+z_{d+1}\alpha_{01+})\log \eta_{t-2}^2+k_{t-3}
% \end{eqnarray*}
% where
%$h_{t-3}, k_{t-3}\in {\cal F}_{t-3}$, from which we deduce that
%$x_3=0.$
Proceeding by induction, we show that ${\bf x}={\bf 0}$ and ${\bf z}={\bf 0}$.
%The arguments used to prove the invertibility of ${\bf J}$ in
%Theorem 4.2 of FWZ (2013) allow to conclude that ${\bf z}={\bf 0}$.
Finally, ${\bf y}={\bf 0}$ and the invertibility of ${\bcalJ}$ is
established.

It follows from Steps {\it i)-v)} and (\ref{4.3}) that
$$\mathbf{S_n}\stackrel{d}{\to} {\cal
N} (\mathbf{0},(\kappa_4 -1){\bcalI}).$$
It can also be shown that
$\widehat{\bcalI} \to {\bcalI}$ and $\widehat{\kappa}_4 \to \kappa_4$ in
probability, from which the conclusion follows. \zak
\section{Test of EGARCH(1,1)}\label{sec5}
In this section, we consider testing the EGARCH(1,1) specification
in the framework of Model (\ref{logGARCHandEGARCH}) with $p=\ell=1$.
For convenience, we reparameterize it as follows
\begin{equation}\label{logGARCHandEGARCH1}
\left\{\begin{array}{lll}
\epsilon_t&=&\sigma_t\eta_t,\quad \\
\log \sigma_t^2&=&\omega_0+\sum_{i=1}^q\omega_{0,i-}1_{\{\epsilon_{t-i}<0\}} + \gamma_0 \eta_{t-1}+\delta_0|\eta_{t-1}|+\beta_{0}\log \sigma_{t-1}^2\\
&&+\sum_{i=1}^q\left(\alpha_{0,i+}1_{\{\epsilon_{t-i}>0\}}+\alpha_{0,i-}1_{\{\epsilon_{t-i}<0\}}\right)\log\epsilon_{t-i}^2.
\end{array}\right.
\end{equation}
 Let
$\bvartheta_0=(\bzeta_0', \balpha_0')'$ where $\bzeta_0=(\omega_0,
\gamma_{0}, \delta_{0}, \beta_{0})'$ and $\balpha_0=(\bomega_{0-}',\balpha_{0+}',
\balpha_{0-}')'$. The vector $\bzeta_0$ is assumed to belong to some
compact parameter set $\Xi\subset \mathbb{R}^4$.

We will derive a LM approach to test the hypothesis that, in
(\ref{logGARCHandEGARCH1}),
$$H_0^{\balpha}: \balpha_0=0\quad\mbox{against}\quad H_1^{\balpha}:
\balpha_0\neq 0.$$ Assuming that  $ |\beta_0|<1$, there exists a
stationary solution to Model (\ref{logGARCHandEGARCH1}) under
$H_0^{\balpha}$, obtained from the MA($\infty$) representation
$$%\begin{equation}\label{eq:Log_var_MA}
\log\sigma_{t}^{2}=\omega_0(1-\beta_0)^{-1}+\sum_{k=1}^{\infty}\beta_0^{k-1}\{\gamma_0 \eta_{t-k}+\delta_0|\eta_{t-k}|\}.
$$%\end{equation}
An important difficulty in the estimation of the EGARCH(1,1) model
is that invertibility
 is not trivial. Invertibility is required to write
$\widetilde{\sigma}_t^2(\bzeta)$, to be defined below, in function
of the observations $\epsilon_t$ for any $\bzeta=(\omega, \gamma,
\delta, \beta)'$. Wintenberger (2013) obtained the following
sufficient condition for continuous invertibility of the
EGARCH(1,1): the compact set $\Xi$ is included in
$\mathbb{R}\times\{\delta\ge |\gamma|\}\times \mathbb{R}^+$ and
$\forall \bzeta\in \Xi$,
\begin{equation} E\left[\log\left(\max \left[\beta , \frac12(\gamma
\epsilon_{0}+\delta |\epsilon_{0}|)\exp\left\{-\frac{\omega}{2(1-\beta
)}\right\}-\beta \right]\right)\right]<0.\label{eq:condtheta}
\end{equation}
Notice that this condition depends on the distribution of the
observations $(\epsilon_t)$.

Denote by $\widehat{\bvartheta}_n^c$ the constrained (by
$H_0^{\balpha}$) estimator of $\bvartheta_0$, defined by
$$%\begin{equation}\label{qmlconst1}
  \widehat{\bvartheta}_n^c=  (\widehat{\bzeta}_n', \mathbf{0}_{1\times 3q})'
%  \argmin_{\bvartheta\in\Theta}\widetilde{Q}_n(\btheta),
$$%\end{equation}
where $\widehat{\bzeta}_n$ is the QMLE of the EGARCH parameters
defined by
$$%\begin{equation}\label{qmlegarch}
  \widehat{\bzeta}_n=\argmin_{\bzeta\in\Xi}\widetilde{Q}_n(\bzeta),
$$%\end{equation}
with
$$
  \widetilde{Q}_n(\bzeta) =n^{-1}\sum_{t=r_0+1}^n\widetilde{\ell}_t(\bzeta),\qquad
 \widetilde{\ell}_t(\bzeta)= \frac{\epsilon_t^2}{\widetilde{\sigma}_t^2(\bzeta)} +\log \widetilde{\sigma}_t^2(\bzeta),
$$
where $r_0$ is a fixed integer and $\log
\widetilde{\sigma}_t^2(\bzeta)$ is recursively defined %, for
%$t=1,2,\dots,n$,
 by
$$\log \widetilde{\sigma}_t^2(\bzeta)=\omega+  \gamma \widetilde{\eta}_{t-1}(\bzeta)+\delta|\widetilde{\eta}_{t-1}(\bzeta)|
+\beta\log \widetilde{\sigma}_{t-1}^2(\bzeta), \quad
\widetilde{\eta}_{t-1}(\bzeta)=\epsilon_{t-1}/\widetilde{\sigma}_{t-1}(\bzeta)$$
using initial values for
$\epsilon_0,\widetilde{\sigma}_{0}^2(\bzeta).$ For any $\bzeta\in \Xi$, the continuous invertibility condition
(\ref{eq:condtheta}) allows to define the sequence $({\sigma}_t^2(\bzeta))_{t\in \mathbb{Z}}$ by
$$\log {\sigma}_t^2(\bzeta)=\omega+  \gamma {\eta}_{t-1}(\bzeta)+\delta|{\eta}_{t-1}(\bzeta)|
+\beta\log {\sigma}_{t-1}^2(\bzeta), \quad
{\eta}_{t-1}(\bzeta)=\epsilon_{t-1}/{\sigma}_{t-1}(\bzeta).$$
We introduce the following assumption.
\begin{itemize}
\item[ {\bf A9:}]
\hspace*{1em} $\bzeta_0 \in \stackrel{\circ}{\Xi}$,
$E(\eta_0^4)<\infty$ and $E\{\beta_0 -\frac12(\gamma_0
\eta_0+\delta_0\left|\eta_0\right|)\}^2<1$.
\end{itemize}
  The following result was established by
Wintenberger (Theorem 6, 2013).
\begin{theo}[{\bf Asymptotics of the QMLE for the EGARCH(1,1)}]\label{qmlEG11}
 For any  compact subset
$\Xi$ of $\mathbb{R}\times\{\delta\ge |\gamma|\}\times \mathbb{R}^+$
satisfying \eqref{eq:condtheta},  almost surely
$\widehat{\bzeta}_n\to \bzeta_0$ as $n\to\infty$ under
$H_0^{\balpha}$. If, in addition, {\bf A9} holds, we have $\sqrt n
(\widehat\bzeta_n-\bzeta_0)\stackrel{d}{\to}\mathcal
N(\bzero,(\kappa_4-1)\bf V^{-1})$ as $n\to\infty$, where ${\bf
V}=E[\nabla \log \sigma_t^2(\bzeta_0) \nabla \log
\sigma_t^2(\bzeta_0)']$ is a positive definite matrix.
\end{theo}
Now, turning to Model (\ref{logGARCHandEGARCH1}), we still denote
 by $\log
\widetilde{\sigma}_t^2(\bvartheta)$ the variable recursively defined,
  for any $\bvartheta$ in $\Xi
\times \mathbb{R}^{3q}$  and $t=1,2,\dots,n$, by
\begin{eqnarray*}\log \widetilde{\sigma}_t^2(\bvartheta)&=&\omega+\sum_{i=1}^q\omega_{i-}1_{\{\epsilon_{t-i}<0\}}+
(\gamma \epsilon_{t-1}+\delta |\epsilon_{t-1}|) e^{-\frac12\log
\widetilde{\sigma}_{t-1}^2(\bvartheta)}
\\
&&+\beta\log
\widetilde{\sigma}_{t-1}^2(\bvartheta)
+\sum_{i=1}^q\left(\alpha_{i+}1_{\{\epsilon_{t-i}>0\}}+\alpha_{i-}1_{\{\epsilon_{t-i}<0\}}\right)\log\epsilon_{t-i}^2,
\end{eqnarray*}
using positive initial values for $\epsilon_0,\dots,
\epsilon_{1-q},\widetilde{\sigma}_{0}^2(\bvartheta)$.

For any $\bvartheta=(\bzeta',  \mathbf{0}_{1\times 3q})',$ the
random vector $\widetilde{\mathbf{D}}_t(\bvartheta)=\frac{\partial}{\partial
\balpha} \log\widetilde{\sigma}_t^2(\bvartheta)$ satisfies
\begin{eqnarray}%\frac{\partial}{\partial \balpha}\log\widetilde{\sigma}_t^2(\bvartheta)
\label{grad1}
\widetilde{\mathbf{D}}_t(\bvartheta)&=&
\widetilde{U}_{t-1}(\bvartheta)\widetilde{\mathbf{D}}_{t-1}(\bvartheta)
%\left(\widetilde{D}_{t-1,1}(\bzeta) ,\ldots,\widetilde{D}_{t-1,2q}(\bzeta) ,{\bf 0}_q\right)'%\frac{\partial}{\partial \balpha} \log\widetilde{\sigma}_{t-1}^2(\bvartheta)
+ \left(
{\bf 1}_{t-1,q}^{-'},\bepsilon_{t-1,q}^{+'},\bepsilon_{t-1,q}^{-'}
\right)'
\end{eqnarray}
where $\widetilde{U}_{t-1}(\bvartheta)=-\frac12  \left\{(\gamma
\epsilon_{t-1}+\delta|\epsilon_{t-1}|)\right\} e^{-\frac12\log
\widetilde{\sigma}_{t-1}^2(\bvartheta)}+\beta$.
%and $\widetilde{D}_{t,i}$ denotes the $i$-th component of $\widetilde{\mathbf{D}}_{t}$, for $i=1, \ldots 3q$,  $t\ge 0$.

Similar to what was accomplished for the Log-GARCH, we will derive
the asymptotic distribution of
\begin{eqnarray*}
\frac1{\sqrt{n}}\sum_{t=1}^n \frac{\partial}{\partial \bvartheta} \widetilde{\ell}_t(\widehat{\bvartheta}_n^c)=
\left(\begin{array}{c}
\mathbf{0}_{4\times 1}\\\mathbf{T}_n:=
\frac1{\sqrt{n}}\sum_{t=1}^n (1-\widehat{\eta}_t^2)\widetilde{\mathbf{D}}_t(  \widehat{\bvartheta}_n^c)
%{\cal B}_{\widehat{\bvartheta}_n^c}^{-1}(B)\widehat{\boldeta}_{t-1,\ell}
%\widetilde{\by}_{t-1,\ell}(\widehat{\bvartheta}_n^c)
\end{array}\right),
\end{eqnarray*}
where
$\widehat{\eta}_t=\epsilon_{t}/\widetilde{\sigma}_t(\widehat{\bvartheta}_n^c).$
Let  $\widehat{\kappa}_4 -1=n^{-1}
\sum_{t=1}^n (1-\widehat{\eta}_t^2)^2,$
\begin{eqnarray*}
\widehat{\bcalK}_{11} &=& \frac1n \sum_{t=1}^n
\widetilde{\mathbf{D}}_t(\widehat{\bvartheta}_n^c)\widetilde{\mathbf{D}}_t(\widehat{\bvartheta}_n^c)'-\left(\frac1n \sum_{t=1}^n
\widetilde{\mathbf{D}}_t(\widehat{\bvartheta}_n^c)\right)\left(\frac1n \sum_{t=1}^n
\widetilde{\mathbf{D}}_t(\widehat{\bvartheta}_n^c)'\right),\\
\widehat{\bf V}&=&\frac1n \sum_{t=1}^n \nabla \log
\widetilde{\sigma}_t^2(\widehat{\bvartheta}_n^c)\nabla' \log
\widetilde{\sigma}_t^2(\widehat{\bvartheta}_n^c),\quad
\widehat{\bPsi}=\frac1n \sum_{t=1}^n \widetilde{\mathbf{D}}_t(\widehat{\bvartheta}_n)\nabla' \log
\widetilde{\sigma}_t^2(\widehat{\bvartheta}_n^c),
\\
\widehat{\bcalK}_{12} &=& -\left\{\widehat{\bPsi}-\left(\frac1n \sum_{t=1}^n
\widetilde{\mathbf{D}}_t(\widehat{\bvartheta}_n)\right)\left(\frac1n \sum_{t=1}^n
\nabla' \log
\widetilde{\sigma}_t^2(\widehat{\bvartheta}_n^c)\right)\right\}\widehat{\bf V}^{-1}=\widehat{\bcalK}_{21}',
\end{eqnarray*}
 and
$$\widehat{\bcalL}=\widehat{\bcalK}_{11}+ \widehat\bPsi \widehat{\bf V}^{-1}\widehat\bPsi'+ \widehat{\bcalK}_{12}\widehat\bPsi'+ \widehat\bPsi \widehat{\bcalK}_{21}.$$
\begin{theo}[{\bf Asymptotic distribution of the LM test under $H_0^{\balpha}$}]\label{test1}
Under the assumptions of Theorem \ref{qmlEG11} (including {\bf A9}),
and under $H_0^{\balpha}$ the matrix $\widehat{\bcalL}$ converges in
probability to a positive definite matrix ${\bcalL}$ and we have
$$\mathbf{LM}_n^{\balpha}=(\widehat{\kappa}_4 -1)^{-1}\mathbf{T}_n'\widehat{\bcalL}^{-1}\mathbf{T}_n\stackrel{d}{\to} \chi^2_{3q}.$$
\end{theo}

\noindent {\bf Proof:} See the supplementary document.\zak

\section{Portmanteau goodness-of-fit tests}
\label{sec6} Portmanteau tests based on residual autocorrelations
are routinely employed in time series analysis, in particular for
testing the adequacy of an estimated ARMA$(p,q)$ model  (see Box and
Pierce (1970), Ljung and Box (1979) and McLeod (1978) for the
pioneer works, and see Li (2004) for a reference book on the
portmanteau tests). The intuition behind these portmanteau tests  is
that if a given time series model with iid innovations $\eta_t$
 is appropriate for the data at hand,
 the autocorrelations of the residuals $\widehat{\eta}_t$ should not be to far from zero.
% The standard portmanteau tests thus consist in rejecting the adequacy of the model for large values of some quadratic form of the residual autocorrelations.

For an ARCH-type model such as Model (\ref{logGARCH}), the portmanteau tests based on residual autocorrelations are irrelevant because we have
$\widehat{\eta}_t=({\sigma_t}/{\widehat{\sigma}_t})\eta_t$ and any process of the form $\epsilon_t=\sigma_t^*\eta_t$, with $\sigma_t^*$ independent of $\sigma\left(\{\eta_u,u< t\}\right)$, is a martingale difference, and thus is uncorrelated. For ARCH-type models, Li and Mak (1994) and Ling and Li (1997) proposed portmanteau tests based on the autocovariances of the {\it squared} residuals.
Berkes, Horv\'ath and Kokoszka (2003)  developed a sharp analysis of the asymptotic theory of these portmanteau tests in the standard GARCH framework (see also Theorem~8.2 in Francq and Zakoïan, 2011). Escanciano (2010) developed diagnostic tests for a general class of conditionally heteroskedastic time series models.
Carbon and Francq (2011) considered the portmanteau tests for the APARCH models.  Recently, Leucht, Kreiss and Neumann (2015) proposed a consistent specification test for GARCH(1,1) models based on a test statistic of
Cramér-von Mises type. The Log-GARCH model is not covered by these works.

To test the null hypothesis
$$H_0: \mbox{the process $(\epsilon_t)$ satisfies Model (\ref{logGARCHSTABLE}),}$$
define the autocovariances of the squared residuals at lag $h$, for $|h|<n$, by
$$\widehat{r}_h=\frac{1}{n}\sum_{t=|h|+1}^n(\widehat{\eta}_t^2-1)(\widehat{\eta}_{t-|h|}^2-1),\qquad \widehat{\eta}_t^2=\frac{\epsilon_t^2}{\widehat{\sigma}^2_t},$$
where  $\widehat{\sigma}_t=\widetilde{\sigma}_t(\widehat{\btheta}_n)$.
For any fixed integer $m$,
$1\leq m<n$,
consider the statistic
$\boldsymbol{\widehat{r}}_m=\left(\widehat{r}_1,\dots,\widehat{r}_m\right)'.$
 Define the $m\times d$ matrix $\widehat{\bK}_m$ whose row $h$, for $1\leq h\leq m$, is the transpose of
\begin{equation}
\label{km}
\widehat{\bK}_m(h,\cdot)=\frac1n
\sum_{t=h+1}^n(\widehat{\eta}_{t-h}^2-1)
\nabla\log\widetilde{\sigma}_t^2(\widehat{\btheta}_n).\end{equation}
The following assumption is marginally milder than {\bf A8}.
\begin{itemize}
\item[{\bf A10:}] The support of $\eta_0$ contains
at least three positive values or three negative values.
\end{itemize}
%\hspace*{1em}
\begin{theo}[{\bf Adequacy test for the AS-Log-GARCH$(p,q)$ model}]
\label{LiMakPortmanteau} Under $H_0$, the assumptions of Theorem
\ref{normality} and {\bf A10}, the matrix $\widehat{\bD}=(\widehat{\kappa}_4-1)^2\bI_m-(\widehat{\kappa}_4-1)\widehat{\bK}_m\widehat{\bJ}^{-1}\widehat{\bK}_m'$ converges in
probability to a positive definite matrix ${\bD}$ and we have
$$n\boldsymbol{\widehat{r}}_m'\widehat{\bD}^{-1}\boldsymbol{\widehat{r}}_m\stackrel{d}{\to} \chi^2_m.$$
\end{theo}
\noindent {\bf Proof:} See the supplementary document.\zak
%Denoting by $\chi_{m}^2(\alpha)$ the $\alpha$-quantile of the chi-square distribution with $m$ degrees of freedom, the adequacy of the Log-GARCH$(p,q)$ model (\ref{logGARCH}) is then rejected at the asymptotic level $\alpha$
%when
%%\begin{equation}
%%\label{portmanteautest}
%$\left\{n\boldsymbol{\widehat{r}}_m'\widehat{\bD}^{-1}\boldsymbol{\widehat{r}}_m>\chi_{m}^2(1-\alpha)\right\}.$
%%\end{equation}
The same result could be established for testing adequacy of an EGARCH(1,1), under {\bf A10} and the assumptions of Theorem
\ref{qmlEG11}.
As usual in portmanteau tests, the choice of $m$ impacts the power of the test. A large $m$ is likely to offer power for a large set of alternatives. Conversely,
choosing $m$ too large may reduce the power for a specific assumption, in particular because the autocovariances will be poorly estimated for large lags.

\section{An application to exchange rates}
\label{sec7}
In the supplementary document,  we investigate the empirical size and power of the LM and portmanteau tests by means of  Monte Carlo simulation experiments.
%We then apply the two tests to assess the relevance of the Log-GARCH(1,1) and EGARCH(1,1)  models on series of exchange rate returns.
We now consider
returns series of the daily exchange rates of the American Dollar
(USD), the Japanese Yen (JPY), the British Pound (BGP), the Swiss
Franc (CHF)  and Canadian Dollar (CAD) with respect to the Euro. The
observations cover the period from January 5, 1999 to January 18,
2012, which corresponds to 3344 observations. The data were obtained
from the web site
\linebreak[1]
\noindent\url{http://www.ecb.int/stats/exchange/eurofxref/html/index.en.html}.

It may seem surprising to investigate asymmetry models for exchange rate returns, while
the conventional view is that leverage is not relevant for such series.
However, many empirical studies (e.g. Harvey and Sucarrat
(2014)), show that asymmetry/leverage is  relevant for exchange rates,
especially when
one currency is
more liquid or more attractive than the other.
It may also be worth mentioning
%that the asymmetry/leverage parameter need not be positive, since
the sign of the
effect depends on which currency appears in the denominator of the exchange rate.

Table~\ref{EstimTaux} displays the
estimated AS-Log-GARCH(1,1) and EGARCH(1,1) models for each series. In order to have two models with the same number of parameters, which facilitates their comparison, we imposed $\alpha=\alpha_{1+}=\alpha_{1-}$ in the AS-Log-GARCH model (see the complementary file for unrestricted estimation of the AS-Log-GARCH(1,1)).
The estimated models are rather similar over the different series. In particular, for the two models and all the series, the
persistence parameter $\beta$ is very high. For all the estimated AS-Log-GARCH models, except the GBP, the value of $\omega_{-}$ is significantly positive, which reflects the existence of a leverage effect.
The leverage effect is also visible in the EGARCH models, because the estimated value of $\gamma$ is negative, except again for the GBP. Comparing the estimated coefficients  $\omega_{-}$ and $\gamma$ with their estimated standard deviations (given in parentheses), the evidence for the presence of a leverage effect is however often weaker in the EGARCH than in the Log-GARCH model. The two models having the same number of parameters, it makes sense to prefer the model with the higher likelihood, given by the last column of Table~\ref{EstimTaux} in bold face. According to this criterion, the Log-GARCH(1,1) is preferred for the USD and GBP series, whereas the EGARCH(1,1) is preferred for the 3 other series.

Even if, for a given series,  a model produces a better fit than the other candidate, this does not guarantee its relevance for that series.
We thus assess the models by means of the two adequacy tests studied in the present paper.
Tables~\ref{AdeqPort} and \ref{AdeqLM}
display the $p$-values of %, respectively,
the portmanteau and LM tests  for testing the null of a AS-Log-GARCH(1,1) (without assuming $\alpha=\alpha_{1+}=\alpha_{1-}$)  and the null of an EGARCH(1,1). The $p$-values smaller than 0.01 are printed in light face. The two tests clearly reject the AS-Log-GARCH(1,1) model for the series JPY, CHF and CAD. The portmanteau tests also clearly reject the EGARCH(1,1) model for the series CHF, and they also find some evidence against the EGARCH(1,1) model for the series JPY, GPD and CAD.  The LM tests  finds strong evidence against the EGARCH(1,1), for all the series except CAD. Using the two adequacy tests, one can  thus arguably reject the EGARCH(1,1) for all the series. Out-of-sample prediction exercises, presented in the supplementary document,  confirm the general superiority of the AS-Log-GARCH over the EGARCH  model for fitting and predicting these series.

To summarize our empirical investigations, the AS-Log-GARCH(1,1) model seems to be relevant for the USD and GBP series, whereas none of the two models is suitable for the 3 other series.

\begin{table}\caption{\label{EstimTaux}
	AS-Log-GARCH(1,1) and EGARCH(1,1) models fitted by QMLE on daily returns of exchange
	rates.}
	\begin{center}
\resizebox{1\textwidth}{!} {			
\begin{tabular}{lccccc}
				\hline\hline\vspace*{-0.3cm}\\
				&\multicolumn{4}{c}{AS-Log-GARCH(1,1)}\\
				\cline{2-6}\vspace*{-0.3cm}\\
Currency &  $\widehat{\omega}$ & $\widehat{\omega}_-$ & $\widehat{\alpha}$ & $\widehat{\beta}$& \mbox{Log-Lik}\\
				USD  &  \phantom{-}0.005      (0.008)    &  \phantom{-}0.037         (0.008)     & 0.021          (0.003)     & 0.972          (0.005)   &     \bf{-0.102}    \\
				JPY  &  \phantom{-}0.022      (0.013)    &  \phantom{-}0.059         (0.013)     & 0.041          (0.005)     & 0.946          (0.007)   &     -0.350     \\
				GBP  &  \phantom{-}0.033      (0.010)    &  -0.003        (0.011)     & 0.030          (0.004)     & 0.964          (0.006)   &     \bf{\phantom{-}0.547}     \\
				CHF  &  -0.025     (0.017)    &  \phantom{-}0.138         (0.017)     & 0.033          (0.005)     & 0.961          (0.006)   &     \phantom{-}1.507     \\
				CAD  &  \phantom{-}0.010      (0.008)    &  \phantom{-}0.021         (0.008)     & 0.020          (0.003)     & 0.971          (0.006)   &     -0.170     \\
				%\cline{2-5}
				\\
				&\multicolumn{4}{c}{EGARCH(1,1)}\\
				\cline{2-6}\vspace*{-0.3cm}\\
				&  $\widehat{\omega}$ & $\widehat{\gamma}$ & $\widehat{\delta}$ & $\widehat{\beta}$& \mbox{Log-Lik}\\
				USD  & -0.119     (0.021)    & -0.017       (0.011)     &  0.131       (0.023)     &  0.981        (0.006)   & -0.100    \\
				JPY  & -0.116     (0.017)    & -0.068       (0.013)     &  0.133       (0.020)     &  0.978        (0.005)   & \bf{-0.333}    \\
				GBP  & -0.306     (0.040)    & \phantom{-}0.004        (0.018)     &  0.289       (0.036)     &  0.945        (0.012)   & \phantom{-}0.529     \\
				CHF  & -0.152     (0.027)    & -0.078       (0.016)     &  0.124       (0.023)     &  0.977        (0.005)   & \bf{\phantom{-}1.582}     \\
				CAD  & -0.079     (0.014)    & -0.007       (0.009)     &  0.089       (0.016)     &  0.988        (0.004)   & \bf{-0.161}    \\
				\hline
			\end{tabular}
}
		\end{center}
	\end{table}

\begin{table}\caption{\label{AdeqPort}
The $p$-values of the
portmanteau adequacy tests.}
\begin{center}
\resizebox{1\textwidth}{!} {			
\begin{tabular}{lcccccccccccc}
\hline\hline\vspace*{-0.3cm}\\
Currency &\multicolumn{12}{c}{$m$}\\
&  1 & 2 & 3 & 4&5&6&7&8&9&10&11&12\\
&\multicolumn{12}{c}{AS-Log-GARCH(1,1)}\\
\cline{2-13}\vspace*{-0.3cm}\\
                           USD  & 0.031 &                     0.095 &                     0.194 &    0.039 &      0.015 &           0.012 &    0.02  &      0.034 &    0.036 &       0.047 &     0.071 &        0.086         \\
\textcolor{myblue}{\sout{JPY}}  & 0.029 & \textcolor{myblue}{0.000} & \textcolor{myblue}{0.000} &    \textcolor{myblue}{0.000} &   \textcolor{myblue}{0.000} & \textcolor{myblue}{0.000} & \textcolor{myblue}{0.000} & \textcolor{myblue}{0.000} &        \textcolor{myblue}{0.000} &         \textcolor{myblue}{0.000} &         \textcolor{myblue}{0.000} &        \textcolor{myblue}{0.000}         \\
 GBP  & 0.020  &    0.014 &    0.012 &    0.012 &   0.017 &       0.033 &  0.041 &         0.064 &    0.077 &         0.111 &     0.121 &       0.143         \\
\textcolor{myblue}{\sout{CHF}}  & \textcolor{myblue}{0.000} & \textcolor{myblue}{0.000} &    \textcolor{myblue}{0.000} &   \textcolor{myblue}{0.000} &  \textcolor{myblue}{0.000} &  \textcolor{myblue}{0.000} &   \textcolor{myblue}{0.000} &    \textcolor{myblue}{0.000} &          \textcolor{myblue}{0.000} &         \textcolor{myblue}{0.000} &         \textcolor{myblue}{0.000} &        \textcolor{myblue}{0.000}         \\
 \textcolor{myblue}{\sout{CAD}}  & \textcolor{myblue}{0.000} &   \textcolor{myblue}{0.000} &  \textcolor{myblue}{0.000} &  \textcolor{myblue}{0.000} &  \textcolor{myblue}{0.000} &   \textcolor{myblue}{0.000} &    \textcolor{myblue}{0.000} &   \textcolor{myblue}{0.000} &          \textcolor{myblue}{0.000} &        \textcolor{myblue}{0.000} &    \textcolor{myblue}{0.000} &    \textcolor{myblue}{0.000}         \\
\\
&\multicolumn{12}{c}{EGARCH(1,1)}\\
\cline{2-13}\vspace*{-0.3cm}\\
%&\multicolumn{7}{c}{$m$}\\
%Currency &  1 & 2 & 3 & 4&6&8&10&12\\
USD  &    0.496 &    0.163 &    0.192 &    0.314 &     0.446 &      0.575 &    0.396 &  0.235 &    0.263 &        0.305 &      0.249 &    0.295      \\
\textcolor{myblue}{\sout{JPY} ?}  &   0.484 &    0.052 &    0.066 &    0.054 &   0.048 &      0.015 &     0.025 &     0.039 &      0.010  &    \textcolor{myblue}{0.002}  &    \textcolor{myblue}{0.002}  &        \textcolor{myblue}{0.004}  \\
\textcolor{myblue}{\sout{GBP} ?}  &   0.195 &    0.013  &    \textcolor{myblue}{0.005}  &\textcolor{myblue}{0.009} & \textcolor{myblue}{0.008} &   \textcolor{myblue}{0.004} & \textcolor{myblue}{0.008} &  \textcolor{myblue}{0.007} &\textcolor{myblue}{0.010} &     0.016 &  0.026 &    0.039    \\
\textcolor{myblue}{\sout{CHF}}   &  \textcolor{myblue}{0.002} &   \textcolor{myblue}{0.003} & \textcolor{myblue}{0.000} & \textcolor{myblue}{0.000} & \textcolor{myblue}{0.000} & \textcolor{myblue}{0.000} & \textcolor{myblue}{0.000} & \textcolor{myblue}{0.000} & \textcolor{myblue}{0.000} & \textcolor{myblue}{0.000} & \textcolor{myblue}{0.000} & \textcolor{myblue}{0.000}     \\
\textcolor{myblue}{\sout{CAD} ?}  &   \textcolor{myblue}{0.006} &    0.020  &    0.050  &    0.089 &    0.094 &    0.121 &    0.114 &      0.126 &     0.179 &       0.241 &      0.313 &      0.390       \\
\hline
\end{tabular}
}
\end{center}
\end{table}

\begin{table}
\begin{center}\caption{\label{AdeqLM}%\footnotesize
The $p$-values of the
LM adequacy tests.  }
\resizebox{1\textwidth}{!} {			
\begin{tabular}{lcccccccccccc}
\hline\hline\vspace*{-0.3cm}\\
Currency &\multicolumn{12}{c}{$\ell$ or $q$}\\
&  1 & 2 & 3 & 4&5&6&7&8&9&10&11&12\\
&\multicolumn{12}{c}{AS-Log-GARCH(1,1)}\\
\cline{2-13}\vspace*{-0.3cm}\\
USD  &    0.895 &    0.951 &    0.818 &    0.932 &     0.884 &     0.852 &     0.877 &     0.831 &     0.865 &     0.864 &     0.599 &     0.589\\
\textcolor{myblue}{\sout{JPY}}  & 0.761 &    0.080  &    \textcolor{myblue}{0.000}     &    \textcolor{myblue}{0.000}     &     \textcolor{myblue}{0.000}     &     \textcolor{myblue}{0.000}     &     \textcolor{myblue}{0.000}     &     \textcolor{myblue}{0.000}     &     \textcolor{myblue}{0.000}     &     \textcolor{myblue}{0.000}     &     \textcolor{myblue}{0.000}    &     \textcolor{myblue}{0.000}         \\
 GBP                            & 0.902 &    0.767 &    0.481 &    0.474 &     0.421 &     0.550  &     0.581 &     0.613 &     0.627 &     0.704 &     0.655 &     0.679     \\
\textcolor{myblue}{\sout{CHF}}  & \textcolor{myblue}{0.000}     &    \textcolor{myblue}{0.000}     &    \textcolor{myblue}{0.000}     &    \textcolor{myblue}{0.000}     &     \textcolor{myblue}{0.000}     &     \textcolor{myblue}{0.000}     &     \textcolor{myblue}{0.000}     &     \textcolor{myblue}{0.000}     &     \textcolor{myblue}{0.000}     &     \textcolor{myblue}{0.000}     &     \textcolor{myblue}{0.000}     &     \textcolor{myblue}{0.000}         \\
\textcolor{myblue}{\sout{CAD}}  & 0.895 &    \textcolor{myblue}{0.004} &    \textcolor{myblue}{0.002} &    \textcolor{myblue}{0.002} &     \textcolor{myblue}{0.002} &     \textcolor{myblue}{0.001} &     \textcolor{myblue}{0.002} &     \textcolor{myblue}{0.005} &     \textcolor{myblue}{0.008} &     0.015 &     0.023 &     0.034     \\
\\
&\multicolumn{12}{c}{EGARCH(1,1)}\\
\cline{2-13}\vspace*{-0.3cm}\\
\textcolor{myblue}{\sout{USD}?}  &   0.461 &    0.067 &    \textcolor{myblue}{0.009} &    0.037 &     0.049 &     0.088 &     0.071 &     0.068 &     0.122 &     0.024 &     \textcolor{myblue}{0.001} &     \textcolor{myblue}{0.000}   \\
\textcolor{myblue}{\sout{JPY}} &   \textcolor{myblue}{0.000} &    \textcolor{myblue}{0.000} &    \textcolor{myblue}{0.000} &    \textcolor{myblue}{0.001} &     \textcolor{myblue}{0.003} &     \textcolor{myblue}{0.004} &     \textcolor{myblue}{0.004} &     \textcolor{myblue}{0.002} &     \textcolor{myblue}{0.004} &     \textcolor{myblue}{0.007} &     \textcolor{myblue}{0.004} &     \textcolor{myblue}{0.002}   \\
\textcolor{myblue}{\sout{GBP}}  &   0.676 &    \textcolor{myblue}{0.006} &    \textcolor{myblue}{0.000} &    \textcolor{myblue}{0.000} &     \textcolor{myblue}{0.000} &     \textcolor{myblue}{0.000} &     \textcolor{myblue}{0.000} &     \textcolor{myblue}{0.000} &     \textcolor{myblue}{0.000} &     \textcolor{myblue}{0.000} &     \textcolor{myblue}{0.000} &     \textcolor{myblue}{0.000}   \\
\textcolor{myblue}{\sout{CHF}}  &   \textcolor{myblue}{0.000} &    \textcolor{myblue}{0.000} &    \textcolor{myblue}{0.000} &    \textcolor{myblue}{0.000} &     \textcolor{myblue}{0.000} &     \textcolor{myblue}{0.000} &     \textcolor{myblue}{0.000} &     \textcolor{myblue}{0.000} &     \textcolor{myblue}{0.000} &     \textcolor{myblue}{0.000} &     \textcolor{myblue}{0.000} &     \textcolor{myblue}{0.000}   \\
 CAD  &   0.112 &    0.128 &    0.031 &    0.034 &     0.037 &     0.059 &     0.119 &     0.203 &     0.308 &     0.400 &     0.469 &     0.409    \\
\hline
\end{tabular}}
\end{center}
\end{table}

\section{Conclusion}
\label{sec8}
The EGARCH and AS-Log-GARCH models do not require any a priori restriction on the parameters because the positivity of the variance is automatically satisfied.
This is often consider as the main advantage of such models, by comparison with other GARCH-type formulations designed to capture the leverage effect.
In empirical applications, the EGARCH model is clearly preferred by the practitioners, the Log-GARCH model being rarely considered. The conclusions of our study are not in accordance with this
predominance.
First, we noted that the two models may produce the same volatility process, though they do not produce the same returns process.
Second, it is now well known that invertibility of the EGARCH requires stringent non explicit conditions.
If such conditions are neglected, results obtained from the statistical inference may be dubious.
Third, the adequacy tests developed in this paper show that the two volatility models are not interchangeable for a given series.
Finally,
our estimation results on real exchange rate data do not allow to validate the EGARCH model for any of the
series under consideration. For the AS-Log-GARCH model, the conclusions are mixed: two over six series passed all adequacy tests, and the
out-of-sample performance is generally superior than that of the EGARCH.

%\bigskip
%\begin{center}
%{\Large\bf References}
%\end{center}

%\end{document}
\clearpage

\appendix
%\end{document}
\begin{center}
\Large{\bf Goodness-of-fit tests for  log and exponential GARCH models: complementary results}
\end{center}
\setcounter{page}{1}
This document contains additional results, in particular illustrations and proofs, that have been removed from the main document to save place.
\section{Illustration to Lemma~\ref{EGARCHisLogGARCH}}

\begin{center}
\begin{figure}%[h]
\vspace*{8.cm} \protect
\includegraphics{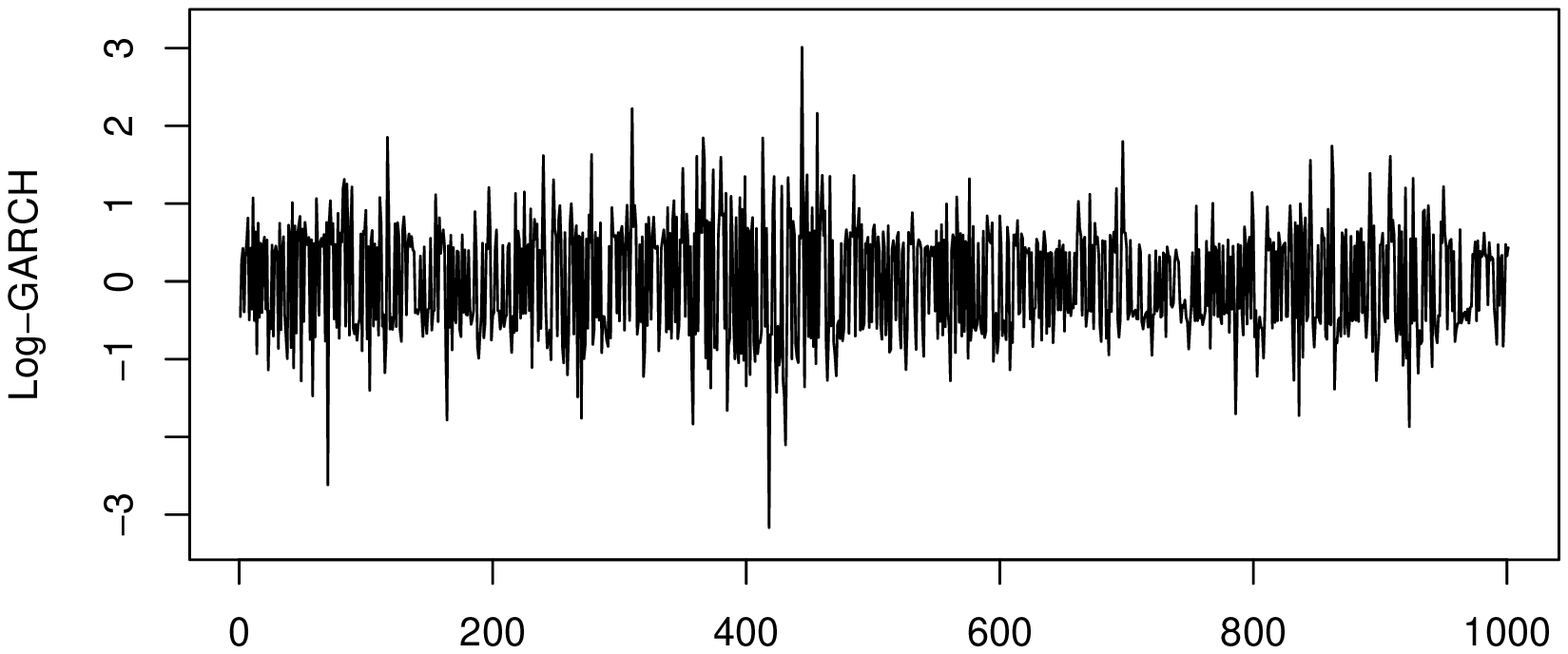}
\vspace*{5.cm}
\includegraphics{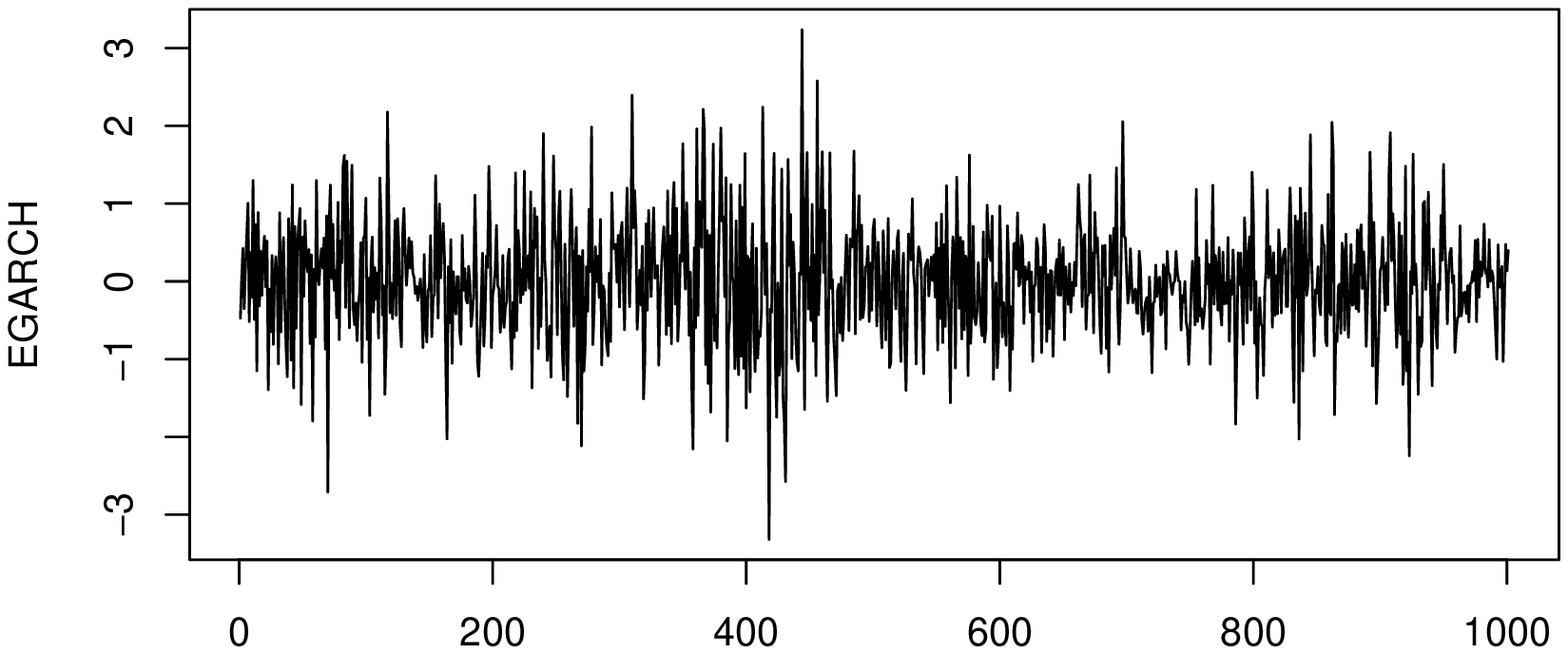}
\vspace*{5.cm}
\includegraphics{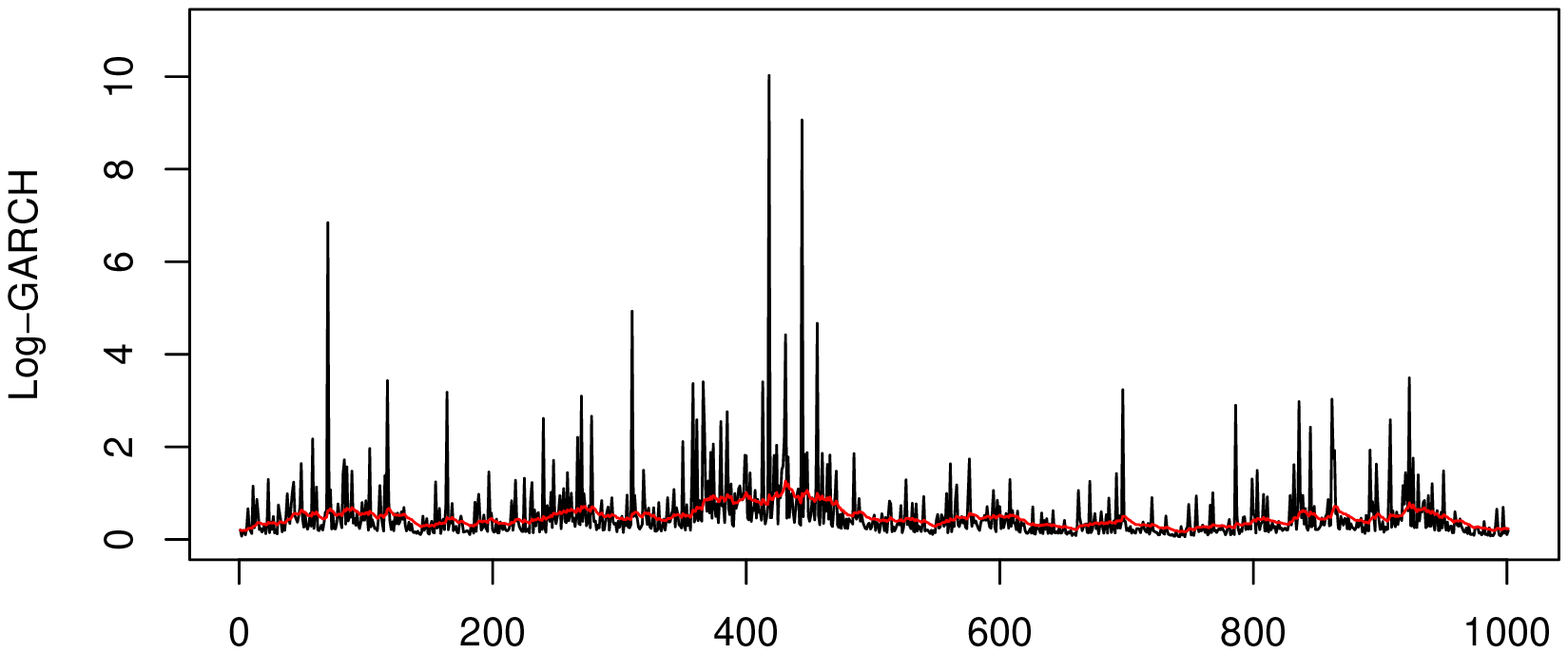}
\vspace*{5.cm}
\includegraphics{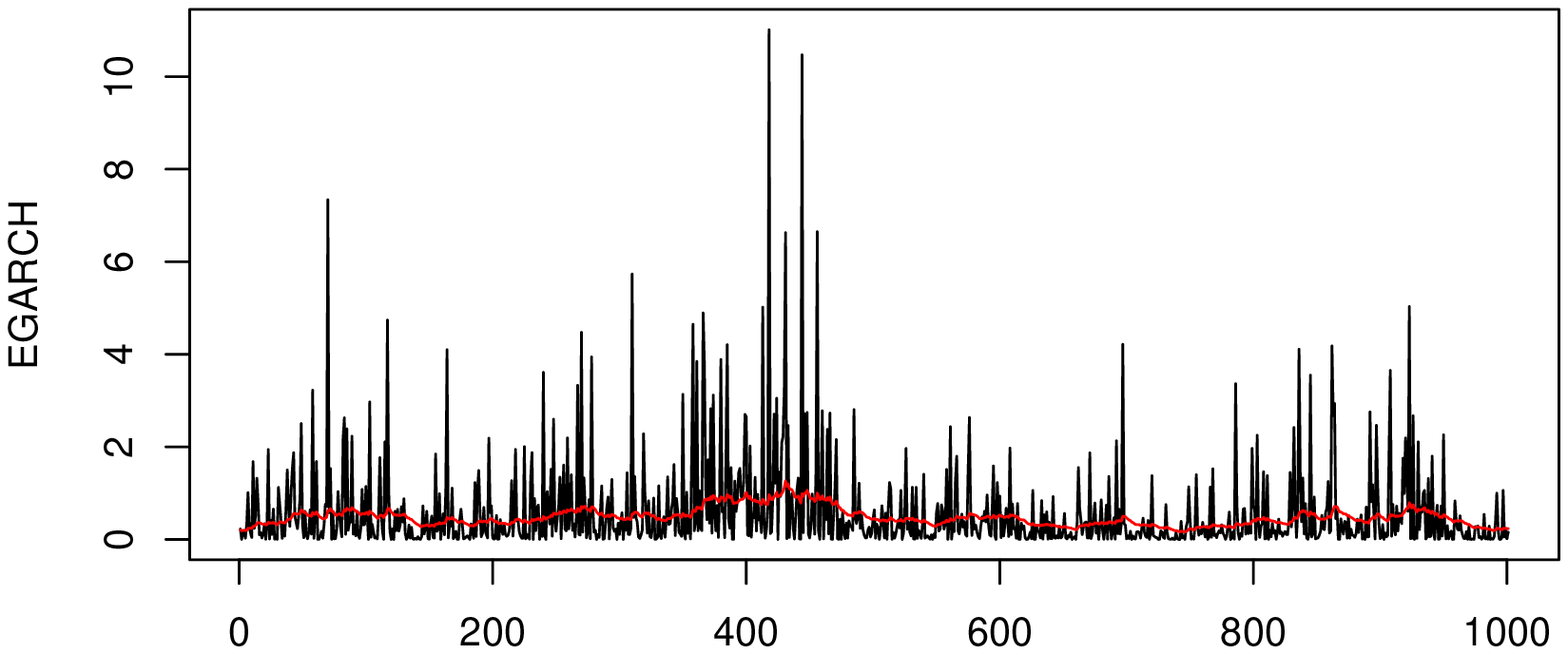}
\vspace{-2cm}
\caption{\label{fig2} {\small Symmetric Log-GARCH(1,1) and EGARCH(1,1) with the same volatility process $\omega=0.2$,
$\alpha=0.2$ and
$\beta=0.95$. The top two panels display the sample paths of the return processes. The
bottom two panels display the sample paths of the squared return processes and the volatilities (in red).}}
\end{figure}
\end{center}

  Note that,  in Lemma~\ref{EGARCHisLogGARCH} for the symmetric case (when $\gamma:=\gamma_+=\gamma_-$), one can take $\alpha=\alpha_+=\alpha_-=\gamma$, $\omega=\tilde{\omega}+\alpha\log Ee^{|\tilde\eta_1|}$, $\beta=\tilde\beta-\gamma$ and
$$\eta_t=\frac{e^{\frac{|\tilde\eta_1|}{2}}}{\sqrt{E e^{|\tilde\eta_1|}}}\mbox{sign}(\tilde\eta_t).$$
  Note also that,  there is a linear relation between $\log(\eta_0^{2})$ and $\tilde \eta_0$ for $\tilde \eta_0\geq 0$, and another linear relation for $\tilde \eta_0< 0$. The tail of  $\eta_t$
is thus heavier than that of $\tilde{\eta}_t$.  This implies that the tails of the Log-GARCH process $\varepsilon_t=\sigma_t\eta_t$ are less impacted by the tails of the volatility process than those of the EGARCH process $\tilde \varepsilon_t=\sigma_t\tilde \eta_t$, leading to possibly  less temporal dependence. To illustrate this point, we plot in Figure~\ref{fig2} trajectories of Log-GARCH(1,1) and EGARCH(1,1) processes with the same symmetric  log volatility process and $\eta_0$ following a standard gaussian distribution. The trajectories have the same periods of high volatilities but the EGARCH(1,1) trajectory looks more blurry when the volatility is low.
%\begin{center}
%\begin{figure}%[h]
%\vspace*{7.cm} \protect
%\special{psfile=lgarchegarch.eps hoffset= 40
%voffset=0 hscale=50 vscale=70}
%\caption{\label{fig2old} {\small Symmetric Log-GARCH(1,1) and EGARCH(1,1) with the same volatility process $\omega=0.2$,
%$\alpha=0.2$ and
%$\beta=0.95$.}}
%\end{figure}
%\end{center}
%\begin{center}
%\begin{figure}%[h]
%\vspace*{7.cm} \protect
%\special{psfile=LGARCH.eps hoffset= 30
%voffset=0 hscale=80 vscale=70}
%\vspace*{8.cm}
%\special{psfile=EGARCH.eps hoffset= 30
%voffset=0 hscale=80 vscale=70}
%\caption{\label{fig2old} {\small Symmetric Log-GARCH(1,1) and EGARCH(1,1) with the same volatility process $\omega=0.2$,
%$\alpha=0.2$ and
%$\beta=0.95$.}}
%\end{figure}
%\end{center}
%The different behaviors in periods of low volatilities could explain the different invertibility properties of the two specifications
%Even if, as discussed before,  the Log-GARCH and EGARCH models have apparent similarities and certain common properties, they are not equivalent.
%In particular, it is obvious to inverse the Log-GARCH model, {\it i.e.} to express the volatility as an explicit function of the past returns,
%whereas the EGARCH(1,1) is invertible only under strong restrictions on the parameters. It is a drawback for the statistical inference of the second specification, see Wintenberger (2013) and FWZ (2013). So one could prefer the first specification, as the Log-GARCH model can produce the same volatility process as the EGARCH (see Lemma~\ref{EGARCHisLogGARCH}) and the inference is much easier with the Log-GARCH. However, the two models are not compatible for a same series and one has to discuss if one specification is more likely to fit the data at hand than the other. It is therefore of interest to develop testing procedures for one specification against the other. This constitutes the main aim of the present paper.

\section{Monte Carlo experiments}

%To assess the performance of the QMLE and
To assess the ability of the adequacy tests to distinguish the two
models, we made the following numerical illustrations.   We
generated $N=1, 000$ independent simulations of length $n=1, 000$ and $n=4, 000$ of a Log-GARCH(1,1) model
with parameter $\btheta_0=(0.01,0.02,0.04,0.05,0.95)$ and an EGARCH(1,1) model
with parameter $\bzeta_0=(-0.15,-0.08,0.12,0.95)$,  both with $\eta_t\sim {\cal
N}(0,1)$. The values of the parameters $\btheta_0$ and $\bzeta_0$ are close to those estimated on the real series of the next section.
On each simulated series, we applied 4 adequacy tests: the LM and portmanteau tests for the null of a Log-GARCH(1,1) and  for the null of an EGARCH(1,1).

Table~\ref{DGPLGARCH} displays
the empirical relative frequencies of rejection over the $N$ replications for the 3 nominal levels $\alpha=1\%$, $5\%$ and $10\%$, when the DGP is the Log-GARCH(1,1) model. Table~\ref{DGPEGARCH} displays the same empirical relative frequencies of rejection when the DGP is the EGARCH(1,1) model.
Recall that, for a random sample of size 1,000, the empirical relative frequency of rejection should vary respectively within the intervals $[0.3;1.9]$, $[3.3;6.9]$ and $[7.6;12.5]$ with probability 0.99 under the assumption that the true probabilities of rejection are respectively $1\%$, $5\%$ and $10\%$.
Tables~\ref{DGPLGARCH} and \ref{DGPEGARCH} show that, as expected the error of first kind is better controlled when $n=4,000$ than when $n=1,000$, both with the LM and portmanteau tests.
The powers of the two tests are quite satisfactory when the null is the Log-GARCH(1,1) model. Even for the sample size $n=1,000$, the two tests are able to clearly reject the Log-GARCH(1,1) model when the DGP is the EGARCH(1,1). For the the null of an EGARCH(1,1), the two tests are less powerful.
For testing the two null assumptions, the LM test is slightly more powerful for small values of $l$ (say $l\leq 4$) whereas the portmanteau test works slightly better with relatively large values of $m$ (say $m\geq 7$).

\begin{table}
\begin{center}\caption{\label{DGPLGARCH}\small
Portmanteau and LM adequacy  tests of the Log-GARCH(1,1) and EGARCH(1,1) models when the DGP is a Log-GARCH(1,1) model.}
%\resizebox{1\textwidth}{!} {
{\small
\begin{tabular}{llcccccccccccc}
\hline\hline\vspace*{-0.5cm}\\
&&\multicolumn{12}{c}{$\ell$ or $q$}\\
  & &1&2&3&4&5&6&7&8&9&10&11&12\\
\multicolumn{14}{l}{Lagrange-Multiplier test for the adequacy of the Log-GARCH(1,1)}\\
%\vspace*{-0.1cm} \\
$n=1000$ &  $\alpha=1\%$    &                       2.2 &   3.0 &   2.6 &   2.7 &   3.0 &   3.2 &   3.6 &   3.6 &   3.7 &   3.3 &   3.2 &   3.2 \\
         &  $\alpha=5\%$    &                       4.8 &   6.4 &   6.3 &   6.5 &   6.5 &   7.2 &   7.4 &   8.0 &   8.3 &   8.8 &   9.0 &   8.7 \\
         &  $\alpha=10\%$   &7.6  &    9.3  &    9.9  &    10.2 &    12.0 &    11.8 &    11.8 &    12.3 &    13.1 &    13.1 &    13.8 &    13.2 \\
\multicolumn{14}{l}{Portmanteau test for the adequacy of the Log-GARCH(1,1)}\\
%\vspace*{-0.1cm} \\
$n=1000$ &  $\alpha=1\%$    &  2.5 &   2.7 &   2.8 &   3.0 &   3.2 &   3.6 &   3.3 &   3.8 &   3.9 &   3.7 &   3.9 &   4.0                          \\
         &  $\alpha=5\%$    &  7.1 &   7.3 &   6.7 &   7.3 &   7.6 &   6.9 &   7.5 &   7.5 &   7.2 &   7.3 &   7.2 &   6.8                          \\
         &  $\alpha=10\%$   &    12.1 &    13.0 &    12.2 &    12.4 &    13.1 &    12.4 &    13.0 &    12.4 &    12.6 &    11.4 &    11.7 &    11.9 \\
\multicolumn{14}{l}{Lagrange-Multiplier test for the adequacy of the EGARCH(1,1)}\\
%\vspace*{-0.1cm} \\
$n=1000$ &  $\alpha=1\%$    &     99.9 &    99.9 &    99.8 &    99.8 &    99.8 &    99.8 &    99.7 &    99.7 &    99.7 &    99.7 &    99.7 &    99.7                        \\
         &  $\alpha=5\%$    &   100 &     100 &     99.9  &     99.9  &     99.8  &     99.9  &     99.9  &     99.9  &     99.9  &     99.8  &     99.8  &     99.8    \\
         &  $\alpha=10\%$   &   100 &     100 &     100 &     100 &     99.9  &     100 &     99.9  &     100 &     100 &     100 &     100 &     100   \\
\multicolumn{14}{l}{Portmanteau test for the adequacy of the EGARCH(1,1)}\\
%\vspace*{-0.1cm} \\
$n=1000$ &  $\alpha=1\%$    &  82.2 &    94.0 &    96.1 &    97.7 &    98.0 &    98.4 &    98.8 &    99.1 &    99.1 &    99.2 &    99.2 &    99.4   \\
         &  $\alpha=5\%$    &  93.4 &    95.9 &    97.5 &    98.2 &    98.4 &    98.6 &    99.1 &    99.2 &    99.3 &    99.4 &    99.4 &    99.4   \\
         &  $\alpha=10\%$   &  96.0 &    97.6 &    98.0 &    98.4 &    98.7 &    98.9 &    99.2 &    99.3 &    99.4 &    99.5 &    99.5 &    99.5   \\
         \\
\multicolumn{14}{l}{Lagrange-Multiplier test for the adequacy of the Log-GARCH(1,1)}\\
%\vspace*{-0.1cm} \\
$n=4000$ &  $\alpha=1\%$    &   1.3 &   1.3 &   1.5 &   2.1 &   2.2 &   2.3 &   2.3 &   2.3 &   2.3 &   2.2 &   2.4 &   2.3                         \\
         &  $\alpha=5\%$    &   2.8 &   4.4 &   4.6 &   4.9 &   5.3 &   5.6 &   6.3 &   6.4 &   6.3 &   5.6 &   6.3 &   6.4                          \\
         &  $\alpha=10\%$   &  4.9  &    6.8  &    8.3  &    8.4  &    8.5  &    9.8  &    10.6 &    11.6 &    11.1 &    11.5 &    11.3 &    10.9    \\
\multicolumn{14}{l}{Portmanteau test for the adequacy of the Log-GARCH(1,1)}\\
%\vspace*{-0.1cm} \\
$n=4000$ &  $\alpha=1\%$    & 2.0 &   1.9 &   2.9 &   2.6 &   2.7 &   3.4 &   3.2 &   3.3 &   3.3 &   3.3 &   3.1 &   3.1                                \\
         &  $\alpha=5\%$    & 5.0 &   5.6 &   6.5 &   6.8 &   7.1 &   6.8 &   7.1 &   7.5 &   6.7 &   7.1 &   7.5 &   7.1                               \\
         &  $\alpha=10\%$   & 10.4 &    10.4 &    10.8 &    11.0 &    11.8 &    11.9 &    11.5 &    12.0 &    11.2 &    11.4 &    11.7 &    11.8     \\
\multicolumn{14}{l}{Lagrange-Multiplier test for the adequacy of the EGARCH(1,1)}\\
%\vspace*{-0.1cm} \\
$n=4000$ &  $\alpha=1\%$    & 100 &     100 &     100 &     100 &     100 &     100 &     100 &     100 &     100 &     100 &     100 &     100         \\
         &  $\alpha=5\%$    & 100 &     100 &     100 &     100 &     100 &     100 &     100 &     100 &     100 &     100 &     100 &     100      \\
         &  $\alpha=10\%$   & 100 &     100 &     100 &     100 &     100 &     100 &     100 &     100 &     100 &     100 &     100 &     100   \\
\multicolumn{14}{l}{Portmanteau test for the adequacy of the EGARCH(1,1)}\\
%\vspace*{-0.1cm} \\
$n=4000$ &  $\alpha=1\%$    & 99.4  &     99.9  &     100 &     100 &     100 &     100 &     100 &     100 &     100 &     100 &     100 &     100    \\
         &  $\alpha=5\%$    & 99.7  &     99.9  &     100 &     100 &     100 &     100 &     100 &     100 &     100 &     100 &     100 &     100     \\
         &  $\alpha=10\%$   & 99.8  &     99.9  &     100 &     100 &     100 &     100 &     100 &     100 &     100 &     100 &     100 &     100     \\

\hline
\end{tabular}%}
}
\end{center}
\end{table}

\begin{table}
	\begin{center}\caption{\label{DGPEGARCH}\small
			As Table \ref{DGPLGARCH}, but when the DGP is an EGARCH(1,1) model.}
		%\resizebox{1\textwidth}{!} {
{\small
		\begin{tabular}{llcccccccccccc}
			\hline\hline\vspace*{-0.3cm}\\
			&&\multicolumn{12}{c}{$q$ or $m$}\\
			& &1&2&3&4&5&6&7&8&9&10&11&12\\
			\multicolumn{14}{l}{Lagrange-Multiplier test for the adequacy of the Log-GARCH(1,1)}\\
			%\vspace*{-0.1cm} \\
			$n=1000$ &  $\alpha=1\%$    & 12.7 &    10.6 &    11.1 &    9.2 &    8.3  &    8.8  &    10.1 &    9.1  & 8.4  &    9.0  &    8.2  &    8.0   \\
			         &  $\alpha=5\%$    & 24.7 &    22.6 &    22.1 &    21.6 &    20.8 &    20.7 &    20.9 &    22.2 & 21.6 &    20.6 &    19.4 &    18.8 \\
			         &  $\alpha=10\%$   & 32.8 &    30.2 &    31.5 &    30.2 &    31.3 &    30.9 &    29.6 &    30.4 & 29.9 &    29.3 &    27.8 &    27.5 \\
			\multicolumn{14}{l}{Portmanteau test for the adequacy of the Log-GARCH(1,1)}\\
			%\vspace*{-0.1cm} \\
			$n=1000$ &  $\alpha=1\%$    & 5.5  &    9.1  &    10.5 &    11.9 &    13.6 &    15.1 &    16.7 &    17.2 & 18.5 &    18.6 &    19.3 &    20.4\\
			         &  $\alpha=5\%$    & 13.7 &    19.3 &    21.7 &    25.1 &    27.5 &    28.2 &    29.8 &    30.8 & 32.0 &    32.0 &    33.7 &    32.8\\
			         &  $\alpha=10\%$   & 20.2 &    26.9 &    32.2 &    33.7 &    36.8 &    38.2 &    40.2 &    42.1 & 40.9 &    40.4 &    40.9 &    42.0\\
			\multicolumn{14}{l}{Lagrange-Multiplier test for the adequacy of the EGARCH(1,1)}\\
			%\vspace*{-0.1cm} \\
			$n=1000$ &  $\alpha=1\%$    &   0.7 &   0.9 &   1.1 &   1.1 & 1.1 &   1.4 &   1.1 &   1.1 &   1.3 &   0.7 &   0.7 &   1.0                  \\
			         &  $\alpha=5\%$    &   3.7 &   4.1 &   5.1 &   4.6 & 5.7 &   5.5 &   5.1 &   5.3 &   5.1 &   4.9 &   4.9 &   5.5                 \\
			         &  $\alpha=10\%$   &  7.2  &    8.7&    8.7&   9.8 &10.9 &  10.7 &  11.6 &  10.3 &  10.2 &  10.1 &  10.5 &    9.6\\
			\multicolumn{14}{l}{Portmanteau test for the adequacy of the EGARCH(1,1)}\\
			%\vspace*{-0.1cm} \\
			$n=1000$ &  $\alpha=1\%$    &    1.2 &   1.3 &   1.5 &   1.6 & 2.6 &   2.4 &   2.4 &   2.6 &   2.6 &   2.9 &   3.3 &   3.4                                             \\
			         &  $\alpha=5\%$    &    6.1 &   6.2 &   7.1 &   6.8 & 7.2 &   7.0 &   8.1 &   8.7 &   8.6 &   7.9 &   8.2 &   8.6                       \\
			         &  $\alpha=10\%$   &   11.1 &    11.9 &    12.6 &    12.6 &    13.4 &    13.6 &    13.8 &    14.6 & 13.9 &    13.9 &    13.7 &    13.6   \\
			\\
			\multicolumn{14}{l}{Lagrange-Multiplier test for the adequacy of the Log-GARCH(1,1)}\\
			%\vspace*{-0.1cm} \\
			$n=4000$ &  $\alpha=1\%$    & 59.4 &    52.9 &    47.3 &    45.9 &    44.5 &    43.9 &    40.2 &    39.9 & 38.9 &    39.0 &    38.2 &    38.3   \\
			         &  $\alpha=5\%$    & 76.8 &    70.8 &    67.1 &    68.1 &    65.3 &    64.0 &    62.6 &    61.5 & 60.9 &    59.5 &    59.9 &    60.0   \\
			         &  $\alpha=10\%$   & 84.3 &    79.6 &    76.6 &    76.0 &    74.2 &    74.5 &    73.6 &    71.8 & 72.0 &    70.7 &    70.7 &    70.0   \\
			\multicolumn{14}{l}{Portmanteau test for the adequacy of the Log-GARCH(1,1)}\\
			%\vspace*{-0.1cm} \\
			$n=4000$ &  $\alpha=1\%$    &  24.0 &    33.6 &    46.5 &    54.6 &    60.1 &    64.2 &    67.5 &    68.4 & 71.0 &    70.8 &    72.1 &    73.6     \\
			         &  $\alpha=5\%$    &  39.8 &    54.8 &    64.0 &    71.7 &    76.1 &    79.3 &    81.0 &    83.2 & 83.3 &    84.2 &    85.3 &    85.7     \\
			         &  $\alpha=10\%$   &  51.1 &    64.6 &    73.8 &    80.1 &    83.1 &    86.6 &    86.6 &    88.2 & 88.6 &    89.7 &    90.8 &    91.0     \\
			\multicolumn{14}{l}{Lagrange-Multiplier test for the adequacy of the EGARCH(1,1)}\\
			%\vspace*{-0.1cm} \\
			$n=4000$ &  $\alpha=1\%$    &  0.6 &   0.7 &   0.8 &   0.3 &   0.5 &   0.2 &   0.9 &   0.4 &   0.7 &   0.8 & 1.0 &   0.8 \\
			         &  $\alpha=5\%$    &  2.2 &   3.6 &   4.0 &   4.0 &   3.6 &   4.1 &   4.8 &   4.3 &   4.7 &   4.4 & 4.5 &   3.7 \\
			         &  $\alpha=10\%$   &  4.5 &   6.4 &   7.5 &   7.6 &   7.5 &   8.0 &   9.2 &   9.2 &   9.3 &   9.4 & 8.8 &   8.6 \\
			\multicolumn{14}{l}{Portmanteau test for the adequacy of the EGARCH(1,1)}\\
			%\vspace*{-0.1cm} \\
			$n=4000$ &  $\alpha=1\%$    &  1.5 &   1.6 &   1.4 &   1.3 &   1.3 &   1.2 &   1.3 &   1.6 &   1.5 &   1.5 & 1.2 &   1.6   \\
			         &  $\alpha=5\%$    &  6.1 &   6.1 &   6.1 &   6.0 &   5.6 &   5.1 &   5.6 &   6.0 &   5.7 &   5.6 & 5.6 &   5.1   \\
			         &  $\alpha=10\%$   & 11.4 &  11.1 &  11.3 &  11.8 &  11.7 &  11.7 &  11.2 &  10.2 &  11.0 &  10.9 &10.6 &  10.2   \\
			\hline
		\end{tabular}%}
}
	\end{center}
\end{table}
\section{Complement to the exchange rates study}
\begin{center}
	\begin{figure}%[h]
		\vspace*{7.cm} \protect
		\includegraphics{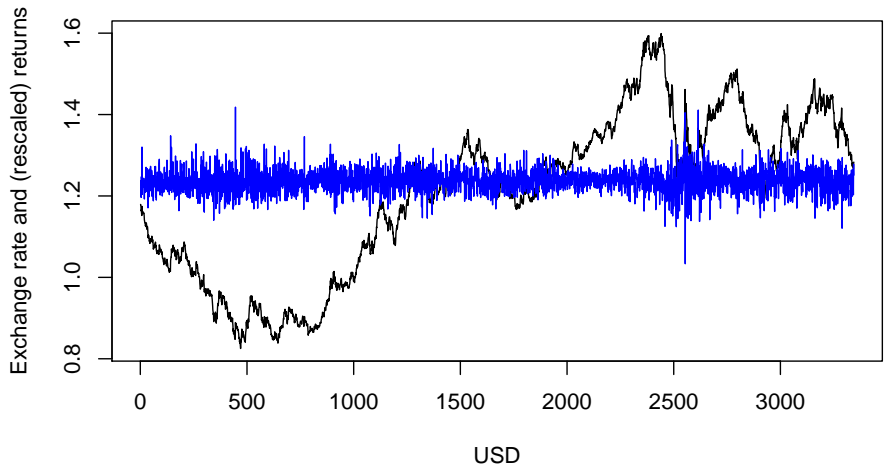}
		\caption{\label{changeUS} {\small Exchange rate and return USD/EURO, from January 5, 1999 to January 18,
				2012.}}
	\end{figure}
\end{center}

Figure~\ref{changeUS} represents the level and return series of the USD to Euro daily exchange rate. %In this figure the returns have been translated
Table~\ref{EstimTauxunrestrict} is the analogue of the top panel of
Table~\ref{EstimTaux}, but for the unrestricted AS-Log-GARCH(1,1).

We also performed out-of-sample predictions of 845 new squared returns, corresponding to the period from January 19, 2012 to May 14, 2015.
As loss function we use either $\left(\epsilon_t^2-\hat{\sigma}_t^2\right)^2$, $\left|\epsilon_t^2-\hat{\sigma}_t^2\right|$, $\left(\log\epsilon_t^2/\hat{\sigma}_t^2\right)^2$, or $\left|\log\epsilon_t^2/\hat{\sigma}_t^2\right|$. Averaging over the 845 observations, we obtain respectively the Mean Squared forecast Errors (MSE), the Mean Absolute forecast Errors (MAE), the MSE of the log-squared returns (log-MSE) and the MAE of the log-squared returns (log-MAE). For the volatility prediction $\hat{\sigma}_t^2$, we used either the Log-GARCH(1,1) or the EGARCH(1,1), both estimated on the initial 3344 observations. Table~\ref{Outsample} shows that the Dielbold-Mariano tests (see Dielbold and Mariano (1995))
often reject the null that the two forecasts are equally accurate in average in favor of the alternative that the EGARCH(1,1) produces less accarate forecasts than the Log-GARCH(1,1), except for the CAD series for which the null can not be rejected.

To summarize our empirical investigations, the Log-GARCH(1,1) model seems to be relevant for the USD and GBP series, whereas none of the two models is suitable for the 3 other series.

%\renewcommand{\baselinestretch}{1}\small\normalsize
%\begin{table}\caption{\label{Outsample}
%$p$-values of the Diebold-Mariano (1995) test for the null that the two models have the same forecast accuracy against the alternative that the Log-GARCH forecasts are less accurate
%or .}
%	\begin{center}
%			\begin{tabular}{lcc cc cc cc}
%\hline
%\\
%				&\multicolumn{2}{c}{MSE}&\multicolumn{2}{c}{MAE}&\multicolumn{2}{c}{Log-MSE}&\multicolumn{2}{c}{Log-MAE}\\
%&Log$\succ$E&E$\succ$Log&Log$\succ$E&E$\succ$Log&Log$\succ$E&E$\succ$Log\\
% USD    &      0.0410 &      0.9590 &      0.0000 &      1.0000 &      0.0000 &      1.0000 &      0.0000 &      1.0000      \\
% JPY    &      0.8131 &      0.1869 &      0.0013 &      0.9987 &      0.0012 &      0.9988 &      0.0064 &      0.9936      \\
% GBP    &      0.6873 &      0.3127 &      0.0322 &      0.9678 &      0.0494 &      0.9506 &      0.1490 &      0.8510      \\
% CHF    &      0.1070 &      0.8930 &      0.0479 &      0.9521 &      0.0000 &      1.0000 &      0.0000 &      1.0000      \\
% CAD    &      0.7373 &      0.2627 &      0.1720 &      0.8280 &      0.0702 &      0.9298 &      0.4524 &      0.5476      \\
%				\hline
%			\end{tabular}
%		\end{center}
%	\end{table}

\begin{table}\caption{\label{EstimTauxunrestrict}
                Unrestricted AS-Log-GARCH(1,1) model fitted by QMLE on daily returns of exchange
                rates.}
                \begin{center}
\resizebox{1\textwidth}{!} {
\begin{tabular}{lcccccc}
                                                               \hline\hline\vspace*{-0.3cm}\\
Currency &  $\widehat{\omega}$ & $\widehat{\omega}_-$ & $\widehat{\alpha}_+$ & $\widehat{\alpha}_-$ & $\widehat{\beta}$& \mbox{Log-Lik}\\
USD  &    \phantom{-}0.008  (0.011)    &    \phantom{-}0.03   (0.011)     &     0.022 (0.005)     &     0.019 (0.005)     &     0.971 (0.005)     &  -0.102 \\
JPY  &    -0.016 (0.017)    &    \phantom{-}0.121  (0.016)     &     0.023 (0.006)     &     0.052 (0.007)     &     0.949 (0.007)     &  -0.343 \\
GBP  &    \phantom{-}0.038  (0.015)    &    -0.014 (0.016)     &     0.031 (0.006)     &     0.027 (0.006)     &     0.965 (0.006)     &  \phantom{-}0.547  \\
CHF  &    -0.135 (0.027)    &    \phantom{-}0.339  (0.027)     &     0.003 (0.006)     &     0.054 (0.007)     &     0.967 (0.005)     &  \phantom{-}1.539  \\
CAD  &    \phantom{-}0.015  (0.011)    &    \phantom{-}0.011  (0.011)     &     0.023 (0.005)     &     0.018 (0.005)     &     0.970 (0.006)     &  -0.170 \\
                                                               \hline
                                               \end{tabular}
}
                               \end{center}
                \end{table}

\begin{table}\caption{\label{Outsample}
$p$-values of the Diebold-Mariano (1995) test for the null that the two models have the same forecast accuracy against the alternative that the EGARCH forecasts are less accurate than those of the Log-GARCH.}
	\begin{center}
			\begin{tabular}{l cc cc c}
\hline
\\
 &USD &    JPY & GBP & CHF & CAD \\
MSE& 0.0410&
0.8131&
0.6873&
0.1070&
0.7373\\
MAE&0.0000 &
0.0013 &
0.0322 &
0.0479 &
0.1720 \\
Log-MSE&0.0000 &
0.0012 &
0.0494 &
0.0000 &
0.0702 \\
Log-MAE&0.0000 &
0.0064 &
0.1490 &
0.0000 &
0.4524 \\
				\hline
			\end{tabular}
		\end{center}
	\end{table}

\section{Proof of Theorem \ref{test1}}
For any $\bvartheta=(\bzeta',  \mathbf{0}_{1\times 3q})',$ let $\mathbf{T}_n({\bzeta})=
\frac1{\sqrt{n}}\sum_{t=1}^n
\{1-{\eta}_t^2({\bzeta})\}\mathbf{D}_t(\bzeta)$ where
${\mathbf{D}}_t(\bzeta)=\frac{\partial}{\partial
\balpha} \log{\sigma}_t^2(\bzeta)$
and let
$\widetilde{\mathbf{T}}_n({\bzeta})= \frac1{\sqrt{n}}\sum_{t=1}^n
\{1-\widetilde{\eta}_t^2({\bzeta})\}\widetilde{\mathbf{D}}_t({\bzeta}), $
where $\eta_t(\bzeta)={\epsilon_t}/{{\sigma}_t(\bzeta)}$ and
$\widetilde{\eta}_t(\bzeta)={\epsilon_t}/{\widetilde{\sigma}_t(\bzeta)}$.

 Define \begin{eqnarray*} {\bcalK}&=&\left(\begin{array}{cc}{\bcalK}_{11}&{\bcalK}_{12}\\{\bcalK}_{21}&{\bcalK}_{22}\end{array}\right),\quad \mbox{where } \quad
{\bcalK}_{11}=\mbox{Var}\{\mathbf{D}_t(\bzeta_0)\},\quad {\bcalK}_{22}=\bf V^{-1}\\ &&\;\;\;\qquad \qquad {\bcalK}_{12}={\bcalK}_{21}'=\mbox{Cov}\{\mathbf{D}_t(\bzeta_0), \nabla \log
\sigma_t^2(\bzeta_0)\} {\bf V}^{-1},  \end{eqnarray*} and
$${\bcalL}={\bcalK}_{11}+ \bPsi{\bf V}^{-1}\bPsi'+ {\bcalK}_{12}\bPsi'+ \bPsi{\bcalK}_{21}$$
where $\bPsi=E\{\mathbf{D}_t(\bzeta_0)\nabla' \log
\sigma_t^2(\bzeta_0)\}.$ Let $\bPsi_i=E\{\mathbf{D}_{t,i}(\bzeta_0)\nabla' \log
\sigma_t^2(\bzeta_0)\}$, where
$\mathbf{D}_{t,i}(\bzeta_0)$ denotes the $i$-th component of
$\mathbf{D}_t(\bzeta_0)$, for $i=1, \ldots 3q$, $t\ge 0$.
Let ${\mathbf{T}}_{n,i}(\bzeta)$
denote the $i$-th component of ${\mathbf{T}}_{n}(\bzeta)$.

The first step of the proof is similar to the one of the proof of Theorem \ref{test}.  Let
$\mathbf{T}_{n,i}$ denote the $i$-th component of
$\mathbf{T}_n=\widetilde{\mathbf{T}}_n(\widehat{\bzeta}_n)$, for $i=1,
\ldots 2\ell$.
 A Taylor
expansion gives, for some ${\bzeta}_*$ between
$\widehat{\bzeta}_n$ and ${\bzeta}_0$,
$$
\mathbf{T}_{n,i} = %=\widetilde{\mathbf{S}}_n(\widehat{\btheta}_n)=
\widetilde{\mathbf{T}}_{n,i} ({\btheta}_0)+\frac1{\sqrt{n}}\frac{\partial \widetilde{\mathbf{T}}_{n,i}}{\partial \bzeta}({\bzeta}_*)
\sqrt{n}(\widehat{\bzeta}_n-
{\bzeta}_0).
$$

We cannot follow the same steps of proof as in Theorem \ref{test} because of the lack of moments in the EGARCH(1,1) model for values of $\bzeta$ satisfying \eqref{eq:condtheta}, see He et al. (2002). However, using the approach of Straumann and Mikosch (2006) refined in Wintenberger (2013), there exist  $K>0$, $\rho\in(0,1)$ and a compact neighborhood ${\cal
V}(\bzeta_0)$ such that
$$
\sup_{\bzeta \in {\cal
V}(\bzeta_0)}|\widetilde{\sigma}_t^2(\bzeta)-\sigma_t^2(\bzeta)|\le K\rho^t,\qquad a.s.
$$
Moreover, the process $\widetilde{\sigma}_t(\bzeta)$  is lower bounded by $\omega/(1-\beta)>0$  under \eqref{eq:condtheta}. By a Lipschitz argument, we then obtain
$$
\sup_{\bzeta \in {\cal
V}(\bzeta_0)}\left|\frac1{\widetilde{\sigma}_t^2(\bzeta)}-\frac1{\sigma_t^2(\bzeta)}\right|\le K\rho^t,\quad a.s.
$$
By an application of Lemma 2.1 in Straumann and Mikosch (2006), it yields to the first assertion (i) below. It remains to
show the three last assertions (ii)-(iv) that are sufficient to prove Theorem \ref{test1}:
\begin{eqnarray*}
&&i) \;  \sup_{\bzeta \in {\cal
V}(\bzeta_0)}\left\|{\mathbf{T}}_n(\bzeta)-
\widetilde{\mathbf{T}}_n(\bzeta)
\right\|\to
0, \quad \sup_{\bzeta \in {\cal
V}(\bzeta_0)}\frac1{\sqrt{n}}\left\|\frac{\partial {\mathbf{T}}_n}{\partial \bzeta}(\bzeta)-
\frac{\partial \widetilde{\mathbf{T}}_n}{\partial \bzeta}(\bzeta)
\right\|\to
0,\\&& \qquad
\;\mbox{ almost surely, }\\
&& ii) \; \left(\begin{array}{c}{\mathbf{T}}_n({\bzeta}_0) \\ \sqrt{n}(\widehat{\bzeta}_n-
{\bzeta}_0)\end{array}\right)\stackrel{d}{\to} {\cal
N} (\mathbf{0},(\kappa_4 -1){\bcalK}), %\quad \mbox{where}
\\&& iii)\;  \frac1{\sqrt{n}} \frac{\partial {\mathbf{T}}_{n,i}}{\partial \bzeta}(\bzeta_{*})\to
\bPsi_i, %=E\{\nabla' \log \sigma_t^2(\btheta_0)\bnu_t(\bvartheta_0^c)\}
\quad
\;\mbox{almost surely, where ${\bzeta}_*$ is between
$\widehat{\bzeta}_n$ and ${\bzeta}_0$,}
 \\
&& v) \; {\bcalL} \;\mbox{ is non-singular.}
%&& vi) \; \widehat{\bcalI} \to {\bcalI} \;\mbox{ in probability.}
\end{eqnarray*}
%To prove i) we note that
%\begin{eqnarray*}
%\left\|{\mathbf{T}}_n(\bzeta)- \widetilde{\mathbf{T}}_n(\bzeta)
%\right\|&\leq &\frac1{\sqrt{n}}
%\sum_{t=1}^n |1-{\eta}_t^2(\bzeta)|\left\|\mathbf{D}_t(\bzeta)-\widetilde{\mathbf{D}}_t(\bzeta)\right\|
%+\frac1{\sqrt{n}}
%\sum_{t=1}^n |\widetilde{\eta}_t^2(\bzeta)-{\eta}_t^2(\bzeta)|\left\|\widetilde{\mathbf{D}}_t(\bzeta)\right\|\\
%&=& T_1(\bzeta)+T_2(\bzeta).
%\end{eqnarray*}
 To prove {\it ii)}, we use that
\begin{eqnarray*}
 \left(\begin{array}{c}{\mathbf{T}}_n({\bzeta}_0) \\ \sqrt{n}(\widehat{\bzeta}_n-
{\bzeta}_0)\end{array}\right)=
\frac1{\sqrt{n}}\sum_{t=1}^n
(1-{\eta}_t^2)\left(\begin{array}{c}\mathbf{D}_t(\bzeta_0)\\
-{\bf V}^{-1}\nabla \log \sigma_t^2(\bzeta_0)
\end{array}\right)+o_P(1).
\end{eqnarray*}
The convergence in distribution thus follows from the central limit
theorem for martingale differences.

The proof of {\it iii)} relies on an almost sure uniform argument applied to $\partial {\mathbf{T}}_n/\partial \bzeta(\bzeta)$ on some neighborhood of $\bzeta_0$. As $\bzeta_\ast$ converges almost surely to $\bzeta_0$, step {\it i)} ensures that
$$
\frac1{\sqrt{n}}\left|\frac{\partial {\mathbf{T}}_n}{\partial \bzeta}(\bzeta_\ast)-
\frac{\partial \widetilde{\mathbf{T}}_n}{\partial \bzeta}(\bzeta_0)
\right|\to
0\qquad a.s.$$ Thus, the result will follow from the ergodic theorem applied to $(\nabla {\bf T}_n(\bzeta_0))$ if ${\bf \Psi}$ is finite. Indeed, the linear stochastic recurrent equation \eqref{grad1} when $\bzeta=\bzeta_0$ takes a simple form with a Lipschitz coefficient equals to $\beta_0-\frac12(\gamma_0\eta_t+\delta_0|\eta_t|)$. Under {\bf A9}, one can use a contractive argument in $L^2$ to prove that $E\{{\bf D}_{t,i}(\bzeta_0)^2\}<\infty$, $i=1,\ldots,3q$. The same argument was already used in Wintenberger (2013) to prove that $E\{\nabla' \log
\sigma_t^2(\bzeta_0)\nabla \log
\sigma_t^2(\bzeta_0)\}<\infty$. Thus, the finiteness of ${\bf \Psi}_i$ is derived from the Cauchy-Schwarz inequality and step {\it iii)} follows.

 Let us prove step {\it iv)}. Suppose there exist ${\bf
x}=(x_i)\in \mathbb{R}^{3q}$ and ${\bf y}\in \mathbb{R}^4$ such that
\begin{equation}\label{eqident1}
{\bf x}'\mathbf{D}_t(\bzeta_0)+{\bf y}' {\bf V}^{-1}\nabla \log
\sigma_t^2(\bzeta_0)=0.\end{equation}
Let ${\bf z}'={\bf y}' {\bf V}^{-1}$.
In view of (\ref{grad1}) we have
\begin{eqnarray*}
\nabla \log
\sigma_t^2(\bzeta_0)&=&{U}_{t-1}(\bzeta_0)\nabla \log
\sigma_{t-1}^2(\bzeta_0)+ \left(1,\epsilon_{t-1}, |\epsilon_{t-1}|,\log \sigma_{t-1}^2(\bzeta_0)\right)',\\
%\left(\begin{array}{c}
%1\\\epsilon_{t-1}\\ |\epsilon_{t-1}|\\\log \sigma_{t-1}^2(\bzeta_0)
%\end{array}\right),\\
{\mathbf{D}}_t(\bzeta_0)&=&
{U}_{t-1}(\bzeta_0){\mathbf{D}}_{t-1}(\bzeta_0)
+\left({\bf 1}_{t-1,q}^{-'},\bepsilon_{t-1,q}^{+'},\bepsilon_{t-1,q}^{-'}\right)'.
\end{eqnarray*}
By stationarity, it follows from (\ref{eqident1}) that
\begin{eqnarray}\label{eq3}
\bx'\left({\bf 1}_{t-1,q}^-,\bepsilon_{t-1,q}^+,\bepsilon_{t-1,q}^-\right)'+
\bz'\left(1,\epsilon_{t-1}, |\epsilon_{t-1}|,\log \sigma_{t-1}^2(\bzeta_0)\right)'=0,\quad a.s.
\end{eqnarray}
It follows that, with notations already used,
\begin{eqnarray}\label{eq2}
x_11_{\{\eta_{t-1}<0\}}+x_{q+1}\log \epsilon_{t-1}^21_{\{\eta_{t-1}<0\}}+
x_{2q+1}\log \epsilon_{t-1}^21_{\{\eta_{t-1}>0\}}&&\nonumber \\+z_1+
z_2\eta_{t-1}\sigma_{t-1}(\bzeta_0)+z_3|\eta_{t-1}|\sigma_{t-1}(\bzeta_0)+
z_4\log\sigma_{t-1}^2(\bzeta_0)&=&R_{t-2}.
\end{eqnarray}
Thus, conditioning on $\eta_{t-1}<0$ we find
\begin{eqnarray*}
x_1+x_{q+1}\log \eta_{t-1}^2
+z_1+
(z_2+z_3)\eta_{t-1}\sigma_{t-1}(\bzeta_0)+
(z_4+x_{q+1})\log\sigma_{t-1}^2(\bzeta_0)&=&R_{t-2}.
\end{eqnarray*}
By arguments already used, in view of Assumption {\bf A8} this entails $x_{q+1}=z_2+z_3=0$.
By conditioning on $\eta_{t-1}>0$ we find $x_{2q+1}=z_2-z_3=0$ and (\ref{eq2}) reduces to
\begin{eqnarray*}
x_11_{\{\eta_{t-1}<0\}}
+z_1+z_4\log\sigma_{t-1}^2(\bzeta_0)&=&R_{t-2}.
\end{eqnarray*}
The sign of $\eta_{t-1}$ being independent of $\sigma\left(\{\eta_u,u\leq t-2\}\right)$ we also have $x_1=0$.
Turning back to (\ref{eq2}), we get %, for some variable $R_{t-3}$ belonging to
%$\sigma\left(\{\eta_u,u\leq t-3\}\right)$,
$
x_21_{\{\eta_{t-2}<0\}}+x_{q+2}\log \epsilon_{t-2}^21_{\{\eta_{t-2}>0\}}+
x_{2q+2}\log \epsilon_{t-2}^21_{\{\eta_{t-2}<0\}}+z_1+
z_4\log\sigma_{t-1}^2(\bzeta_0)=R_{t-3}.
$
Because
$\log\sigma_{t-1}^2(\bzeta_0)=\omega_0+ \gamma_0 \eta_{t-2}+\delta_0|\eta_{t-2}|+\beta_{0}\log \sigma_{t-2}^2(\bzeta_0)$ we get, for
$\eta_{t-2}<0$,
\begin{eqnarray*}
x_2+
x_{2q+2}\log \eta_{t-2}^2+z_1+
z_4(\omega_0+ (\gamma_0-\delta_0) \eta_{t-2})&=&R_{t-3}^*.
\end{eqnarray*}
By arguments already used, we deduce that $x_{2q+2}=z_4=0$.
By conditioning on $\eta_{t-2}>0$, we get $x_{q+2}=0$ and thus $x_2=0$. Proceeding similarly we show that all the components of $\bx$ are equal to zero.
Using (\ref{eq3}), we thus have $z_1=0$. We have shown that, in (\ref{eqident1}), $\bx=\bzero$ and $\by=\bzero$ which entails that
${\bcalL}$ is non-singular.
\zak

\section{Proof of Theorem \ref{LiMakPortmanteau}} Introduce the vector
$\boldsymbol{r}_m=\left(r_1,\dots,r_m\right)'$
where
$$r_h=n^{-1}\sum_{t=h+1}^ns_ts_{t-h},\qquad \mbox{ with  } s_t=\eta_t^2-1\mbox{ and  }0<h<n.$$
Let  $s_t(\btheta)$ (respectively $\widetilde{s}_t(\btheta)$)  be the random variable
obtained by replacing $\eta_t$ by $\eta_t(\btheta)=\epsilon_t/\sigma_t(\btheta)$ (respectively $\widetilde{\eta}_t(\btheta)=\epsilon_t/\widetilde{\sigma}_t(\btheta)$) in $s_t$. Let  $r_h(\btheta)$ (respectively $\widetilde{r}_h(\btheta)$)  be
obtained by replacing $\eta_t$ by $\eta_t(\btheta)$ (respectively $\widetilde{\eta}_t(\btheta)$) in $r_h$. The vectors $\boldsymbol{r}_m(\btheta)=\left(r_1(\btheta),\dots,r_m(\btheta)\right)'$ and $\boldsymbol{\widetilde{r}}_m(\btheta)=\left(\boldsymbol{\widetilde{r}}_1(\btheta),\dots,\boldsymbol{\widetilde{r}}_m(\btheta)\right)'$ are such that $\boldsymbol{r}_m=\boldsymbol{r}_m(\btheta_0)$ and $\boldsymbol{\widehat{r}}_m=\boldsymbol{\widetilde{r}}_m(\widehat{\btheta}_n)$.

We first study the asymptotic impact of the unknown initial values on the statistic $\boldsymbol{\widehat{r}}_m$.
We have $s_t(\btheta)s_{t-h}(\btheta)-\widetilde{s}_t(\btheta)\widetilde{s}_{t-h}(\btheta)=a_t+b_t$ with $a_t=\left\{s_t(\btheta)-\widetilde{s}_t(\btheta)\right\}s_{t-h}(\btheta)$ and $b_t=\widetilde{s}_t(\btheta)\left\{s_{t-h}(\btheta)-\widetilde{s}_{t-h}(\btheta)\right\}$.
A straightforward adaptation of the proof of (\ref{eq1}) shows that the right-hand side can be replaced by $K\rho^t$ in this inequality.
Thus, we have
\begin{eqnarray*}
\left|a_t\right|&\leq &
K\rho^t\epsilon_t^2\left(\frac{\sigma_{t-h}^2}{\sigma_{t-h}^2(\btheta)}\eta_{t-h}^2+1\right).\end{eqnarray*}
Lemma~\ref{lem:mom} and the $c_r$ and Hölder inequalities  entail that for sufficiently small $s^*\in(0,1)$, there exists a neighborhood $\mathcal V $ of $\btheta_0$ such that
$$E\left|\frac{1}{\sqrt{n}}\sum_{t=1}^n\sup_{\btheta\in\mathcal V}|a_t|\right|^{s^*}\leq Kn^{-s^*/2}\sum_{t=1}^n\rho^{ts^*}\to 0$$
as $n\to\infty$. It follows that $n^{-1/2}\sum_{t=1}^n\sup_{\btheta\in\mathcal V}|a_t|=o_P(1)$.
The same convergence holds for $b_t$ and for the derivatives of $a_t$ and $b_t$. We then obtain
\begin{equation}
\label{valinit4port}
\sqrt{n}\left\|\boldsymbol{r}_m-\widetilde{r}_m(\btheta_0)\right\|=o_P(1),\qquad
\sup_{\btheta\in\mathcal V}\left\|\nabla \boldsymbol{r}'_m(\btheta)-\nabla \boldsymbol{\widetilde{r}}'_m(\btheta)\right\|=o_P(1).
\end{equation}

We now show that the asymptotic distribution of $\sqrt{n}\boldsymbol{\widehat{r}}_m$ is a function of the joint asymptotic distribution of
$\sqrt{n}\boldsymbol{r}_m$ and of the QMLE.
Using (\ref{valinit4port}) and the consistency of $\widehat{\btheta}_n$, Taylor expansions of the components of
$\boldsymbol{r}_m(\cdot)$
around
$\widehat{\btheta}_n$
and
$\btheta_0$
shows that
\begin{eqnarray*}
\sqrt{n}\boldsymbol{\widehat{r}}_m&=&\sqrt{n}\boldsymbol{\widetilde{r}}_m(\btheta_0)+\left[\nabla {\widetilde{\boldsymbol{r}}}'_m(\btheta^*)\right]'\sqrt{n}(\widehat{\btheta}_n-\btheta_0)\\
&=&\sqrt{n}\boldsymbol{r}_m+\left[\nabla {\boldsymbol{r}}'_m(\btheta^*)\right]'\sqrt{n}(\widehat{\btheta}_n-\btheta_0)+ o_P(1)
\end{eqnarray*}
where the $h$-th row of the matrix $\left[\nabla {\widetilde{\boldsymbol{r}}}'_m(\btheta^*)\right]'$ is the transpose of $\nabla \widetilde{r}_h(\btheta_h^*)$ for some $\btheta_h^*$ between $\widehat{\btheta}_n$ and $\btheta_0$. In Section 7.11 of FWZ, we have shown the existence of moments of all order for $\log \sigma_t^2(\btheta)$ and their derivatives at any order, uniformly in $\btheta\in{\cal V}$ for some
  neighborhood ${\cal V}$ of $\btheta_0$.
Together with Lemma~\ref{lem:mom}, this implies
  that
$$
%\label{der24port}
E \sup_{\btheta\in{\cal V}}\left|\frac{\partial^2 s_t(\btheta)s_{t-h}(\btheta)}{\partial\btheta_i\partial\btheta_j}\right|<\infty\quad\mbox{for all }i,j\in\{1,\dots,d\}.
$$
Using these inequalities,   the assumption $E\eta_t^4<\infty$,  and the almost sure convergence of $\btheta_h^*$ to $\btheta_0$,  Taylor expansions and the ergodic theorem yield
$$ \nabla r_h(\btheta_h^*)=\nabla r_h(\btheta_0)+o_P(1)\to
c_{h}:=E\left\{s_{t-h}\nabla s_t(\btheta_0)\right\}=-E\left\{s_{t-h}
\nabla\log\sigma_t^2(\btheta_0)\right\}.
$$
Note that $c_{h}$ is the almost sure limit of (\ref{km}). Let $\bK_m$ be the $m\times d$ matrix whose $h$-th row is $c_{h}'$.
%For the next to last equality, we use the fact that $E\left\{s_{t}{\partial s_{t-h}(\btheta_0)}/{\partial
%\btheta}\right\}=0$.
We have shown  that
\begin{eqnarray}
\label{4port}
\sqrt{n}\boldsymbol{\widehat{r}}_m
=\sqrt{n}\boldsymbol{r}_m+\bK_m\sqrt{n}(\widehat{\btheta}_n-\btheta_0)+o_P(1).
\end{eqnarray}

We now derive  the asymptotic distribution of
$\sqrt{n}(\boldsymbol{r}_m,\widehat{\btheta}_n-\btheta_0)$.
Note that
$$\boldsymbol{r}_m=\frac{1}{n}\sum_{t=1}^ns_t\boldsymbol{s}_{t-1:t-m}+o_P(1)\quad\mbox{where}\quad\boldsymbol{s}_{t-1:t-m}=(s_{t-1},\dots,s_{t-m})'.$$
 With this notation, we have $\bK_m=-E\boldsymbol{s}_{t-1:t-m}\nabla' \log\sigma_t^2(\btheta_0)$.
We have seen in the proof of Theorem \ref{test} that
$$
\sqrt{n}\left(\widehat{\btheta}_n-\btheta_0\right)=-\bJ^{-1}\frac{1}{\sqrt{n}}\sum_{t=1}^n (1-\eta^2_t)\nabla \log\sigma_t^2(\btheta_0)+o_P(1).
$$
The central limit theorem applied to the martingale difference
$$\left\{\left(s_t\nabla'\log\sigma_t^2(\btheta_0),s_t\boldsymbol{s}_{t-1:t-m}'\right)';\sigma\left(\eta_u,\, u\leq t\right)\right\}$$
then shows that
\begin{eqnarray}
\label{ANACFQLME}
\sqrt{n}\left(\begin{array}{c}\widehat{\btheta}_n-\btheta_0
\\
\boldsymbol{r}_m\end{array}\right)&=&\frac{1}{\sqrt{n}}\sum_{t=1}^ns_t\left(\begin{array}{c}\bJ^{-1}\nabla \log\sigma_t^2(\btheta_0)\\
\boldsymbol{s}_{t-1:t-m}\end{array}\right)+o_P(1)\nonumber\\
&\stackrel{\bcalL}{\to}&{\cal
N}\left\{\bzero,\left(\begin{array}{cc}(\kappa_4-1)\bJ^{-1}&\bSigma_{\widehat{\btheta}_n\boldsymbol{r}_m}\\\bSigma_{\widehat{\btheta}_n\boldsymbol{r}_m}'&(\kappa_4-1)^2\bI_m\end{array}\right)\right\},
\end{eqnarray}
where
\begin{eqnarray*}
\bSigma_{\widehat{\btheta}_n\boldsymbol{r}_m}&=&
(\kappa_4-1)\bJ^{-1}E\nabla \log\sigma_t^2(\btheta_0)\boldsymbol{s}_{t-1:t-m}'=-(\kappa_4-1)\bJ^{-1}\bK_m'.
\end{eqnarray*}

Using together (\ref{4port}) and  (\ref{ANACFQLME}), we obtain
\begin{eqnarray*}
\sqrt{n}\boldsymbol{\widehat{r}}_m
&\stackrel{\bcalL}{\to}&{\cal
N}\left(\bzero,\bD\right),\quad \bD=(\kappa_4-1)^2\bI_m-(\kappa_4-1)\bK_m\bJ^{-1}\bK_m'.
\end{eqnarray*}

We now show that $\bD$ is invertible. Assumption {\bf A3} entails that the law of $\eta_t^2$ is non degenerated. We thus have $\kappa_4>1$, and it remains to show the invertibility of
$$(\kappa_4-1)\bI_m-\bK_m\bJ^{-1}\bK_m'=E \mathbf{V}\mathbf{V}',\quad \mathbf{V}=\boldsymbol{s}_{-1:-m}+\bK_m\bJ^{-1}\nabla \log\sigma_0^2(\btheta_0).$$ If this matrix were singular then there would exist $\blambda=(\lambda_1,\dots,\lambda_m)'$ such that  $\blambda\neq \bzero$ and
\begin{equation}
\label{Vinversible}
\blambda'\mathbf{V}=\blambda'\boldsymbol{s}_{-1:-m}+\bmu'\nabla \log\sigma_0^2(\btheta_0)=\bzero\quad\mbox{ a.s.,}
\end{equation}
 with $\bmu'=\blambda'\bK_m\bJ^{-1}$. %Note that $\bmu=(\mu_1,\dots,\mu_{d})'\neq \bzero$. Otherwise $\blambda'\boldsymbol{s}_{-1:-m}=\bzero$ a.s., which
%would contradict the fact that  the  $s_t$'s are independent and non degenerated (by Assumption {\bf A3}).
Note that
\begin{equation}\label{eq:nabla}
\nabla \log\sigma_t^2(\btheta)=\sum_{j=1}^p\beta_j \nabla \log\sigma^2_{t-j}(\btheta)+
\left(1, {\bf 1}_{t-1,q}^-\\\bepsilon_{t-1,q}^+,
\bepsilon_{t-1,q}^-,\bsigma^2_{t-1,p}(\btheta)\right)',
\end{equation}
%Denoting by $R_t$ \emph{any}
%random variable which is measurable with respect to
%$\sigma\{\eta_u,\:u\leq t\}$,
Equation~(\ref{Vinversible}) gives
\begin{equation}
\label{equation15}
\blambda'\mathbf{V}=\lambda_1\eta^2_{-1}+\mu_2 1_{\eta_{-1}<0}+\mu_{2+q} 1_{\eta_{-1}>0}\log\epsilon^2_{-1}+\mu_{2+2q} 1_{\eta_{-1}<0}\log\epsilon^2_{-1}+R_{-2}.
\end{equation}
 Thus (\ref{Vinversible}) entails
the two equations
\begin{equation}
\label{equation1}
1_{\eta_{-1}>0}\left\{\lambda_1\eta^2_{-1}+\mu_{2+q} \log\eta^2_{-1}+R_{-2}\right\}=0\quad\mbox{ a.s.}
\end{equation} and
\begin{equation}
\label{equation2}
1_{\eta_{-1}<0}\left\{\lambda_1\eta^2_{-1}+\mu_{2+2q} \log\eta^2_{-1}+R_{-2}\right\}=0\quad\mbox{ a.s.}
\end{equation}
Note that an equation of the form $a x^{2}+b\log|x|+c=0$ cannot have
more than 2 positive roots or more than 2 negative roots, except if
$a=b=c=0$.  By Assumption {\bf A10}, Equations (\ref{equation1}) and
(\ref{equation2}) thus imply $\lambda_1=0$. We thus also have
$\mu_{2+q}=\mu_{2+2q}=0$ and it follows from (\ref{equation15}) that $\mu_2=0$. Given that $\lambda_1=\mu_2=\mu_{2+q}=\mu_{2+2q}=0$,
(\ref{Vinversible}) and (\ref{eq:nabla}) now give
\begin{eqnarray}
\label{equation16}
\nonumber
\blambda'\mathbf{V}&=&\lambda_2\eta^2_{-2}+\mu_31_{\eta_{-2}<0}+\mu_{3+q}1_{\eta_{-2}>0}\log\epsilon^2_{-2}+\mu_{3+2q}1_{\eta_{-2}<0}\log\epsilon^2_{-2}\\
&&+
\mu_{3+3q}\log\sigma^2_{-1}+R_{-3}=0.
\end{eqnarray}
Since
\begin{eqnarray*}\log\sigma^2_{-1}&=&\omega + \omega_{1-}R_{-3}1_{\eta_{-2}<0}+\alpha_{1+}1_{\eta_{-2}>0}(\log\eta^2_{-2}+R_{-3})\\
&&+\alpha_{1-}1_{\eta_{-2}<0}(\log\eta^2_{-2}+R_{-3})+R_{-3},\end{eqnarray*} we have
the two equations
%\begin{equation}
%\label{equation1bis}
$$
1_{\eta_{-2}>0}\left\{\lambda_2\eta^2_{-2}+(\mu_{3+q} +\mu_{3+3q}\alpha_{1+})\log\eta^2_{-2}+R_{-3}\right\}=0\quad\mbox{ a.s.}
$$
%\end{equation}
and
%\begin{equation}
%\label{equation2bis}
$$
1_{\eta_{-2}<0}\left\{\lambda_2\eta^2_{-2}+(\mu_{3+2q} +\mu_{3+3q}\alpha_{1-})\log\eta^2_{-2}+R_{-3}\right\}=0\quad\mbox{ a.s.}
$$
%\end{equation}
By Assumption {\bf A10}, we obtain $$\lambda_2=\mu_{3+q} +\mu_{3+3q}\alpha_{1+}=\mu_{3+2q} +\mu_{3+3q}\alpha_{1-}=0.$$
In view of (\ref{equation16}), it follows that $\mu_3=0$.
By iterating the previous arguments, it can be shown that
%$\lambda_3=0$. Continuing in this way, we show that
$\lambda_1=\cdots=\lambda_m=0$  which leads to a contradiction. The non-singularity
of $\bD$ follows. The proof of the convergence
$\widehat{\bD}\to \bD$ in probability (and even almost surely) as
$n\to\infty$ is omitted. \zak

\end{document}